\documentclass[12pt]{article}
\usepackage[all]{xy}
\usepackage{latexsym, amsfonts, amscd, amssymb, verbatim, amsmath,
enumerate, amsxtra}
\usepackage{graphicx}
\usepackage{indentfirst}
\usepackage[mathscr]{eucal}

\input{tcilatex}
\begin{document}

\title{Toward Homological Structure Theory of Semimodules: On Semirings All
of Whose Cyclic Semimodules Are Projective}
\author{\textbf{\ }S. N. Il'in$^{1}$, Y. Katsov$^{2}$, T.G. Nam$^{3}$ \\
%EndAName
{\footnotesize Sergey.Ilyin@kpfu.ru; katsov@hanover.edu; tgnam@math.ac.vn;}\\
$^{1}${\footnotesize Lobachevsky Institute of Mathematics and Mechanics}\\
{\footnotesize Kazan (Volga Region) Federal University, 18, Kremlevskaja
St., 420008 Kazan, }\\
{\footnotesize Tatarstan, Russia}\\
$^{2}${\footnotesize Department of Mathematics}\\
{\footnotesize Hanover College, Hanover, IN 47243--0890, USA}\\
$^{3}${\footnotesize Institute of Mathematics, VAST}\\
{\footnotesize 18 Hoang Quoc Viet, Cau Giay, Hanoi, Vietnam}}
\date{}
\maketitle

\begin{abstract}
In this paper, we introduce homological structure theory of semirings and
CP-semirings---semirings all of whose cyclic semimodules are projective. We
completely describe semisimple, Gelfand, subtractive, and anti-bounded,
CP-semirings. We give complete characterizations of congruence-simple
subtractive and congruence-simple anti-bounded CP-semirings, which solve two
earlier open problems for these classes of semirings. We also study in
detail the properties of semimodules over Boolean algebras whose
endomorphism semirings are CP-semirings; and, as a consequence of this
result, we give a complete description of ideal-simple CP-semirings.

\textbf{2000} \textbf{Mathematics Subject Classifications}: 16Y60, 16D99,
06A12; 18A40, 18G05, 20M18.


------------------------------------------------------------------------------------

The third author is supported by Vietnam National Foundation for Science and
Technology Development (NAFOSTED) and the Vietnam Institute for Advanced
Study in Mathematics (VIASM).
\end{abstract}

%\setlength{\baselineskip}{16pt}

%{{\bf Key words:} simple semiring, ideal-simple semiring, congruence-simple semiring,
%complete lattice, finitary-complete semirings.}\\
%{{\bf 2000 Mathematics Subject Classification}. 16Y60, 08A30, 18G05,
%18G99.}
\newtheorem{theo}[subsection]{\bf Theorem}
\newtheorem{exa}[subsection]{\bf
Example} \newtheorem{defi}[subsection]{Definition} %
\newtheorem{prop}[subsection]{Proposition} \newtheorem{cor}[subsection]{\bf
Corollary} \newtheorem{lem}{Lemma} %\theoremstyle{definition}

\section{Introduction}

In the modern homological theory of modules over rings, the results
characterizing rings by properties of modules and/or suitable categories of
modules over them are of great importance and sustained interest (for a good
number of such results one may consult, for example, \cite{lam:lomar}).
Inspired by this, quite a few results related to this genre have been
obtained in different nonadditive settings during the last three decades.
Just to mention some of these settings, we note that a very valuable
collection of numerous interesting results\ on characterizations of monoids
by properties and/or by categories of acts over them, \textit{i.e.}, on
so-called \textit{homological classification of monoids}, can be found in a
recent handbook \cite{kilp-kn-mik:maac}; and, for the results on `\textit{%
homological classification of distributive lattices},'\ one may consult the
survey \cite{fof:podl}.

Nowadays, on the other hand, one may clearly notice a growing interest in
developing the algebraic theory of semirings and their numerous connections
with, and applications in, different branches of mathematics, computer
science, quantum physics, and many other areas of science (see, for example,
\cite{golan:sata} and \cite{glazek:agttlos}). As algebraic objects,
semirings certainly are the most natural generalization of such (at first
glance different) algebraic systems as rings and bounded distributive
lattices. Thus, when investigating semirings and their representations one
should undoubtedly use methods and techniques of\ both ring and lattice
theories as well as diverse techniques and methods of categorical and
universal algebra. Thus, the wide variety of the algebraic techniques
involved in studying semirings, and their representations/semimodules,
perhaps explains why the research on homological
characterization/classification of semirings is still behind that for rings
and monoids (for some recent results on `\textit{homological
characterization of semirings},'\ one may consult \cite{asw:cosbpaps}, \cite%
{il'in:svwasai}, \cite{il'in:otatsottfttoram}, \cite{il'in:dsoisadpopsos},
\cite{il'in:v-s}, \cite{kat:tpaieosoars}, \cite{kat:ofsos}, \cite{kat:thcos}%
, \cite{knt:ossss}, \cite{kn:meahcos}, \cite{knt:mosssparp}, \cite{wjk:eras}%
, and \cite{aikn:ovsasaowcsai}).

In the current and previous \cite{aikn:ovsasaowcsai} our papers, we have
studied semirings all of whose cyclic semimodules are projective and
injective, so-called \textit{CP-} and \textit{CI-semirings}, respectively.
Both these papers, of course, belong to the homological characterization of
semirings, but, taking into consideration the additive Yoneda lemma \cite[%
Proposition 1.3.7]{bor:hoca} (see also \cite[Theorem 3.7.1]{macl:cwm}), the
cyclic semimodules permit us to look at the homological characterization
from a little bit different perspectives and, to the extend of our
knowledge, at the very first time explicitly initiate a new, obviously
\textquotedblleft bloody\textquotedblright\ connected with the homological
characterization, approach/area --- `homological structure theory (of
(semi)rings in our case, for instance)'--- (for more detail, please, see
Sections 3 and 6). Perhaps this phenomenon can, at least implicitly, explain
a persistent interest in studying of cyclic modules over rings (see, \textit{%
e.g.}, the recent monograph \cite{jst:cmatsor}).

Also, `congruence- and ideal-simple semirings' constitutes another booming
area in semiring research which has quite interesting and promising
applications in various fields, in particular in cryptography \cite%
{mmr:pkcbosa} (for some relatively recent developments in this area we refer
our potential readers to \cite{bshhurtjankepka:scs}, \cite{bashkepka:css},
\cite{jezkepka:tesoas}, \cite{jezkepkamaroti:tesoas}, \cite{kake:anofgics},
\cite{knt:mosssparp}, \cite{knz:ososacs}, \cite{kz:fsais}, \cite{mf:ccs},
\cite{monico:ofcss}, and \cite{zumbr:cofcsswz}). In this respect, in the present
paper we consider, in the context of CP-semirings, congruence- and
ideal-simple semirings as well. The paper is organized as follows.

For the reader's convenience, all subsequently necessary notions and facts
on semirings and semimodules over them are included in Section 2.

In Section 3, we briefly present some \textquotedblleft
ideological\textquotedblright\ background of the homological structure
theory of (semi)rings (Theorems 3.1, 3.3, and Remark 3.4), establish some
important general properties of CP-semirings (Proposition 3.10) as well as
the \textquotedblleft structural\textquotedblright\ theorem for CP-semirings
(Theorem 3.11), in particular showing that the class of CP-semirings is
essentially wider than that of semisimple rings.

In Section 4, among the main results of the paper we single out the
following ones: Full descriptions of semisimple CP-semirings (Theorem 4.4),
of Gelfand \cite[p. 56]{golan:sata} CP-semirings (Theorem 4.9), subtractive
CP-semirings (Theorem 4.10), and anti-bounded CP-semirings (Theorem 4.16);
Full descriptions of congruence-simple subtractive and anti-bounded
CP-semirings (Corollary 4.11 and Theorem 4.17, respectively).

Main results of Section 5 are the following ones: Full descriptions of
semimodules over Boolean algebras whose endomorphism semirings are
CP-semirings (Theorems 5.9 and 5.10); Full descriptions of ideal-simple
CP-semirings (Theorem 5.11 and Corollary 5.12).

In conclusive Section 6, we briefly presented, in our view, some interesting
and important directions for furthering considerations of this and previous
\cite{aikn:ovsasaowcsai} our papers as well as posted a few specific
problems.

Finally, all notions and facts of categorical algebra, used here without any
comments, can be found in \cite{macl:cwm}; for notions and facts from
semiring and lattice theories, we refer to \cite{golan:sata}, and \cite%
{birkhoff:lathe} (or \cite{gratzer:glt} and \cite{skor:eolt}), respectively.{%
\ {\ }}

\section{Preliminaries}

Recall \cite{golan:sata} that a \textit{semiring}\emph{\/} is an algebra $%
(S,+,\cdot ,0,1)$ such that the following conditions are satisfied:\medskip

(1) $(S,+,0)$ is a commutative monoid with identity element $0$;

(2) $(S,\cdot, 1)$ is a monoid with identity element $1$;

(3) Multiplication distributes over addition from either side;

(4) $0s=0=s0$ for all $s\in S$.\medskip

As usual, a \textit{right\/ }$S$\textit{-semimodule} over the semiring $S$
is a commutative monoid $(M,+,0_{M})$ together with a scalar multiplication $%
(m,s)\mapsto ms$ from $M\times S$ to $M$ which satisfies the identities $%
m(ss^{^{\prime }})=(ms)s^{^{\prime }}$, $(m+m^{^{\prime }})s=ms+m^{\prime }s$%
, $m(s+s^{^{\prime }})=ms+mr^{^{\prime }}$, $m1=m$, $0_{M}s=0_{M}=m0$ for
all $s,s^{^{\prime }}\in S$ and $m,m^{^{\prime }}\in M$.\medskip

\textit{Left semimodules}\emph{\/} over $S$ and homomorphisms between
semimodules are defined in the standard manner. And, from now on, let $%
\mathcal{M}$ be the variety of commutative monoids, and $\mathcal{M}_{S}$
and $_{S}\mathcal{M}$ denote the categories of right and left semimodules,
respectively, over a semiring $S$.

Recall (\cite[Definition 3.1]{kat:tpaieosoars}) that the tensor product
bifunctor $-\otimes -:$ $\mathcal{M}_{S}$ $\times $ $_{S}\mathcal{M}%
\longrightarrow \mathcal{M}$ on a right semimodule $A\in |\mathcal{M}_{S}|$
and a left semimodule $B\in |_{S}\mathcal{M}|$ can be described as the
factor monoid $F/\sigma $ of the free monoid $F\in |\mathcal{M}|$, generated
by the Cartesian product $A\times B$, factorized with respect to the
congruence $\sigma $ on $F$ generated by ordered pair having the form
\begin{equation*}
\langle (a_{1}+a_{2},b),(a_{1},b)+(a_{2},b)\rangle ,\langle
(a,b_{1}+b_{2}),(a,b_{1})+(a,b_{2})\rangle ,\langle (as,b),(a,sb)\rangle ,
\end{equation*}%
with $a_{1},a_{2}\in A$, $b_{1},b_{2}\in B$ and $s\in S$.

For a right $S$-semimodule $M$, we will use the following subsemimodules:%
\begin{equation*}
\begin{tabular}{lll}
$I^{+}(M)$ & $:=$ & $\{m\in M\,|\,m+m=m\};$ \\
$Z(M)$ & $:=$ & $\{z\in M\mid z+m=m\text{ for some }m\in M\};$ \\
$V(M)$ & $:=$ & $\{m\in M\mid m+m^{\prime }=0\text{ for some }m^{\prime }\in
M\}.$%
\end{tabular}%
\end{equation*}
In the special case when $M =S$ (viewed as a right semimodule over itself),
these subsets are (two-side) ideals of $S$. Also, we denote by $%
I^{\times}(S) $ the set of all multiplicatively-idempotent elements of $S$.

A right $S$-semimodule $M$ is \textit{zeroic} (\textit{zerosumfree,
additively idempotent}) if $Z(M)=M$ ($V(M)=0,$ $I^{+}(M)=M$). In particular,
a semiring $S$ is \textit{zeroic} (\textit{zerosumfree, additively
idempotent)} if $S_{S}$ $\in |\mathcal{M}_{S}|$ is a zeroic (zerosumfree,
additively idempotent) semimodule.

A subsemimodule $K$ of a right $S$-semimodule $M$ is \textit{subtractive}%
\emph{\ } if for all $m,m^{\prime }\in M$, $m$ and $m+m^{\prime }\in K$
imply $m^{\prime }\in K$. A right $S$-semimodule $M$ is \textit{subtractive}
if it has only subtractive subsemimodules. A semiring $S$ is \textit{right
subtractive} if $S$ is a subtractive right semimodule over itself. \emph{\ }

As usual (see, for example, \cite[Ch. 17]{golan:sata}), if $S$ is a
semiring, then in the category $\mathcal{M}_{S}$, a \textit{free} (right)
semimodule $\sum_{i\in I}S_{i},S_{i}\cong $ $S_{S}$, $i\in I$, with a basis
set $I$ is a direct sum (a coproduct) of $I$-th copies of $S_{S}$. And a
\textit{projective} right semimodule in $\mathcal{M}_{S}$ is just a retract
of a free right semimodule. A semimodule $M_{S}$ is \textit{finitely
generated} iff it is a homomorphic image of a free semimodule with a finite
basis set. The semimodule $M_{S}$ is \textit{cyclic} iff it is a homomorphic
image of the free semimodule $S_{S}$.

\textit{Congruences} on a right $S$-semimodule $M$ are defined in the
standard manner, and $\mathrm{Cong}(M_{S})$ (or, simply by $\mathrm{Cong}(M)$
when $S$ is fixed) denotes the set of all congruences on $M_{S}$. This set
is non-empty since it always contains at least two congruences---the \textit{%
diagonal congruence} $\vartriangle _{M}:= \{(m,m)$ $|$ $m\in M \}$ and
the \textit{universal congruence}\emph{\ }$M^{2}:=\{(m,n)$ $|$ $m,n\in M$ $\}
$. Any subsemimodule $L$ of a right $S$-semimodule $M$ induces a congruence $%
\equiv _{L}$ on $M$, known as the \textit{Bourne congruence}, by setting $%
m\equiv _{L}m^{\prime }$ iff $m+l=m^{\prime }+l^{\prime }$ for some $%
l,l^{\prime }\in L$.

Following \cite{bashkepka:css}, a semiring $S$ is \textit{congruence-simple}
if the diagonal, $\vartriangle _{S}$, and the universal, $S^{2}$,
congruences are the only congruences on $S$; and $S$ is \textit{ideal-simple}
if $0$ and $S$ are its only ideals. A semiring $S$ is said to be \textit{%
simple} if it is simultaneously congruence-simple and ideal-simple. Note
that in a semiring setting, these notions are not the same (see, \textit{e.g.%
}, \cite[Examples 3.8]{knz:ososacs}) and should be differed.

\section{On homological structure theory of semimodules and CP-semirings}

Let us consider a category $\mathcal{M}_{S}$ of the right $S$-semimodules
over a semiring $S$. Any semimodule $M\in |\mathcal{M}_{S}|$ can be
naturally considered as an additive functor $M:S\longrightarrow \mathcal{M}$
from\ the one object category $S$ to the category of commutative monoids $%
\mathcal{M}$; and thanks mainly to the additive Yoneda lemma \cite[%
Proposition 1.3.7]{bor:hoca}, the only \textquotedblleft representable
functor\textquotedblright\ $\mathcal{M}_{S}(S_{S},-):\mathcal{M}%
_{S}\longrightarrow \mathcal{M}_{S}$, corresponding to the regular right
semimodule $S_{S}\in |\mathcal{M}_{S}|$, produces the familiar natural
(functorial) isomorphism $\mathcal{M}_{S}(S_{S},M)\approx M$ of right $S$%
-semimodules. By using this observation, we obtain the following fundamental
result of the so-called \textquotedblleft \textit{homological structure
theory}\textquotedblright\ of semimodules --- representation of semimodules
as colimits of diagrams of the regular semimodule $S_{S}$ --- which is an
one object additive analog of the colimits of representable functors for
set-valued functors (\textit{cf} \cite[Theorem 3.7.1]{macl:cwm} and \cite[%
Proposition 1.3.8]{bor:hoca}), namely\medskip

\noindent \textbf{Theorem 3.1} \textit{Any semimodule }$M\in |M_{S}|$\textit{%
\ can be represented (in a canonical way) as a colimit of a functor }$D:%
\mathcal{C}\longrightarrow \mathcal{M}_{S}$\textit{\ from a small category }$%
\mathcal{C}$\textit{\ to }$\mathcal{M}_{S}$ \textit{that has the regular
semimodules }$S_{S}$\textit{\ as its values on the objects; in short, a
semimodule }$M\in |\mathcal{M}_{S}|$\textit{\ can be represented (in a
canonical way) as a colimit of a diagram of the regular semimodules }$S_{S}$%
\textit{.}\medskip \textit{\ }

\noindent \textbf{Proof. }Given a semimodule $M\in |\mathcal{M}_{S}|$, to
construct the needed diagram we first consider the \textquotedblleft
category of elements,\textquotedblright\ $\mathrm{Elts}(M)$, of $M$ with
objects $m\ $for each element$\mathbb{\ }m\in M$ and with morphisms $%
s:m\longrightarrow n$ for those elements $s\in S$ for which $ms=n$. Then,
let $\mathcal{C}:=$ $\mathrm{Elts}(M)^{op}$ and $D:\mathcal{C}%
\longrightarrow \mathcal{M}_{S}$ be the functor which sends each object $%
m\in |\mathcal{C}|$ to the regular semimodule $S_{S}^{m}:=S_{S}$ and each
morphism $s\in \mathcal{C}(n,m)$ to the induced homomorphism $s^{\ast
}:S_{S}^{n}\longrightarrow S_{S}^{m}$. Then, using the natural functor
isomorphism $\mathcal{M}_{S}(S_{S},-)\simeq Id_{\mathcal{M}_{S}}$ and
actually word by word repeating the proofs of \cite[Theorem 3.7.1]{macl:cwm}
or \cite[Proposition 1.3.8]{bor:hoca} for our \textquotedblleft one object
additive\textquotedblright\ case, we obtain that $M=Co\lim D$, \textit{i.e.}%
, the semimodule $M$ is a colimit of the (small) diagram of the regular
semimodules $\{S_{S}^{m},S_{S}^{n}\in |\mathcal{M}_{S}|;s^{\ast
}:S_{S}^{n}\longrightarrow S_{S}^{m}$ $|$ $m,n\in M,s\in S$ and $ms=n\}$.%
\textit{\ \ \ \ \ \ }$_{\square }\medskip $

Taking into consideration that, for any $m\in M$, the image of the
semimodule homomorphism $S_{S}\longrightarrow M$, given by $1_{S}\longmapsto
m$, is the cyclic subsemimodule $mS\subseteq M_{S}$ of the semimodule $M_{S}$%
, one can readily rephrase the previous result as\medskip

\noindent \textbf{Corollary 3.2} \textit{A semimodule }$M\in |\mathcal{M}%
_{S}|$\ \textit{can be represented (in a canonical way) as a colimit of a
diagram of all its cyclic subsemimodules, i.e., of the diagram }$\{%
mS,nS\subseteq M_{S};s^{\ast }:nS\longrightarrow mS$ $|$ $m,n\in M,s\in S$
\textit{and the semimodule homomorphism }$s^{\ast }$\textit{is defined by
the equation }$ms=n \}$\textit{.\ \ \ \ \ \ }$_{\square
}\medskip $\textit{\ }

This observation very well explains a strong interest in serious studying of
cyclic (semi)modules (for a good collection of very interesting results on
cyclic modules we refer a reader to the recent monograph \cite{jst:cmatsor}%
), as well as motivates necessity of a developing of the so-called
homological structure theory of semimodules which in the matter of fact is a
studying of colimits of diagrams of cyclic semimodules possessing some
special important properties, and which we explicitly at the first time
initiate in this paper. It is absolutely clear that the homological
structure theory of semimodules is very closely connected with the, nowadays
widely recognized, homological characterization, or classification, of
semirings. Thereby, for example, calling the \textquotedblleft
canonical\textquotedblright\ diagram of all cyclic subsemimodules of a
semimodule $M_{S}$ the \textit{full c-diagram} of $M_{S}$ and the \textit{%
full injective (projective) c-diagram} of $M_{S}$ provided that all objects
of the full c-diagram of $M_{S}$ are injective (projective) semimodules, the
celebrated characterization of semisimple rings given by B. Osofsky (see,
\textit{e.g.}, \cite[Theorem 1.2.9]{lam:afcinr}, or \cite[Corollary 6.47]%
{lam:lomar})\ and \cite[Theorem 1.2.8]{lam:afcinr} can be rephrased as
\medskip

\noindent \textbf{Theorem 3.3} \textit{The following conditions for a ring }$%
S$ \textit{are equivalent:}

\textit{(1) The full c-diagram of every module} $M\in |\mathcal{M}_{S}|$
\textit{is injective;}

\textit{(2) A colimit of the full c-diagram of every module} $M\in |\mathcal{%
M}_{S}|$ \textit{is an injective\ module;}

\textit{(3)} \textit{The full c-diagram of every module} $M\in |\mathcal{M}%
_{S}|$ \textit{is projective;}

\textit{(4)} \textit{A colimit of the full c-diagram of every module} $M\in |%
\mathcal{M}_{S}|$ \textit{is a projective\ module;}

\textit{(3) }$S$\textit{\ is a classical semisimple ring.\medskip }

\noindent \textbf{Proof. }One only should note that every cyclic
(semi)module is always an object of, and a colimit of, its full c-diagram.%
\textit{\ \ \ \ \ \ }$_{\square }\medskip $\textit{\ }

\noindent \textbf{Remark 3.4} However, for semirings in general, the
conditions in Theorem 3.3 are not equivalent: Thus, for the semirings $%
Ext(R) $, introduced in \cite{aikn:ovsasaowcsai}, if $R$ is a classical
semisimple ring, by \cite[Theorem 4.18]{aikn:ovsasaowcsai}, (1) of Theorem
3.3 is true, but, by Corollary 3.2 and \cite[Theorem 3.4]%
{il'in:dsoisadpopsos}, (2) is not true; Also, for a finite Boolean algebra $%
B $ with $|B|$ $>$\textit{\ }$1$, by \cite[Theorem 4.3]{aikn:ovsasaowcsai},
(1) is true, but, by \cite[Theorem 3]{fof:tpopodl} (see also \cite[Section 4]%
{fof:podl}), (2) is not satisfied. Even these observations, for example,
obviously open a wide avenue for furthering, from several different
perspectives, the homological structure theory of semimodules as well as its
\textquotedblleft bloody\textquotedblright\ connection with the homological
characterization of semirings.\textit{\medskip }

As was shown in \cite[Theorem 3.4]{il'in:dsoisadpopsos} and \cite[Theorem 4]%
{il'in:svwasai}, respectively, a semiring $S$ is a (classical) semisimple
ring iff all right $S$-semimodules are projective, iff all finitely
generated right $S$-semimodules are projective. In light of this and taking
into consideration Theorem 3.3, it is quite natural to consider semirings $S$
with the full projective c-diagrams for every semimodule $M\in |\mathcal{M}%
_{S}|$--- in other words, semirings $S$ over which all cyclic right (left)
semimodules are projective. Thus, in this section, we exactly initiate a
study of such semirings, that we call the right (left)\textit{\ CP-semirings}%
, and show that the class of CP-semirings is significantly wider than that
of semisimple rings. First, the following important observation will prove
to be useful.\medskip

\noindent \textbf{Proposition 3.5} \textit{A homomorphic image of a right
(left) CP-semiring is a right (left) CP-semiring itself.\medskip }

\noindent \textbf{Proof. }Let $S$ be a right CP-semiring and $\pi
:S\longrightarrow T$ a surjective homomorphism of semirings. Then, by \cite[%
Section 4]{kat:thcos}, the surjection $\pi $ gives the rise to two functors:
the \textit{restriction} functor $\pi ^{\#}:\mathcal{M}_{T}\longrightarrow
\mathcal{M}_{S}$, defined by $b.s=a\pi (s)$ for any $b\in B\in $ $|\mathcal{M%
}_{T}|$ and $s\in S$; and the \textit{extension} functor $\pi
_{\#}:=-\otimes _{S}\pi ^{\#}T=-\otimes _{S}T:\mathcal{M}_{S}\longrightarrow
\mathcal{M}_{T}$. Obviously, the restriction functor $\pi ^{\#}:\mathcal{M}%
_{T}\longrightarrow \mathcal{M}_{S}$ preserves cyclic semimodules. Then,
using the natural adjunction $\pi _{\#}:\mathcal{M}_{S}\rightleftarrows
\mathcal{M}_{T}:\pi ^{\#}$ (\cite[Proposition 4.1]{kat:thcos}) and the
natural isomorphism of the functors $\pi _{\#}\pi ^{\#}, $ $Id_{\mathcal{M}%
_{T}}:\mathcal{M}_{T}\longrightarrow \mathcal{M}_{T} $ (\cite[Proposition 4.6%
]{kat:thcos}), we get that $T$ is a right CP-semiring.\textit{\ \ \ \ \ \ }$%
_{\square }\medskip $\textit{\ }

From this proposition right away follows\medskip\

\noindent \textbf{Corollary 3.6 }\textit{A semiring }$S=S_{1}\oplus S_{2}$%
\textit{\ is a right (left) CP-semiring iff the semirings }$S_{1}$\textit{\
and }$S_{2}$\textit{\ are right (left) CP-semiring.\medskip }

The following observation is almost obvious.\medskip

\noindent \textbf{Lemma 3.7 }\textit{If }$S$\textit{\ is a right
CP-semiring, then, for every congruence }$\tau $\textit{\ on the right }$S$%
\textit{-semimodule }$S_{S}\in $\textit{\ }$|\mathcal{M}_{S}|$\textit{,
there exist an idempotent }$e\in I^{\times }(S)$\textit{\ and an }$S$\textit{%
-isomorphism }$\psi \colon \overline{S}:=S/\tau \rightarrow eS$\textit{\
such that }$1\tau e$\textit{, }$\psi (\overline{1})=e$\textit{, and }$p\psi
=1_{\overline{S}}$\textit{, where }$p\colon S\rightarrow \overline{S}$%
\textit{\ is the canonical projection.\medskip }

\noindent \textbf{Proof. }Since the right $S$-semimodule $\overline{S}$ is
projective, there exists an $S$-monomorphism $\psi \colon {\overline{S}}%
\rightarrow S$ such that $p\psi =1_{\overline{S}}$, and let $e=\psi (%
\overline{1})$. Then, we have $\overline{1}=p\psi (\overline{1})=p(e)=%
\overline{e}$, that is, $1\tau e$. Hence, $e=\psi (\overline{1})=\psi (%
\overline{e})=\psi (\overline{1})e=e^{2}$ and, as $\psi $ is injective, $%
\overline{S}\cong \psi (\overline{S})=eS$.\textit{\ \ \ \ \ \ }$_{\square
}\medskip $\textit{\ }

\noindent \textbf{Lemma 3.8 }\textit{If }$\equiv _{I}$ \textit{is the Bourne
congruence on a right CP-semiring }$S$\textit{\ defined by a right ideal }$%
I_{S}\subseteq S$\textit{, then there exists an idempotent }$e\in I^{\times
}(S)$\textit{\ such that }$1\equiv _{I}e$\textit{\ and }$I\subseteq $ Ann$%
_{r}(e)$\textit{.\medskip }

\noindent \textbf{Proof. }We may assume that the Bourne congruence $\equiv
_{I}$ is not the universal one --- otherwise, the statement becomes trivial
by taking $e=0$. Then, $\overline{S}:=S/\equiv _{I}$ is a nonzero semimodule
and, by Lemma~3.7, there exists an $S$-isomorphism $\psi \colon \overline{S}%
\rightarrow eS$ such that $\psi (\overline{1})=e$ and $1\equiv _{I}e$.
Clearly, for any $x\in I$, $\overline{1}x=\overline{x}=\overline{0}$, hence,
$ex=\psi (\overline{1})x=\psi (\overline{1}x)=\psi (\overline{0})=0$, that
is, $I\subseteq \mathrm{Ann}_{r}(e)$.\textit{\ \ \ \ \ \ }$_{\square
}\medskip $\textit{\ }

\noindent \textbf{Lemma 3.9 }\textit{Let }$S$\textit{\ be a right
CP-semiring and }$\{e_{i},\,i\in \mathbb{N}\}$\textit{\ a subset of }$%
I^{\times }(S)$\textit{\ such that }$e_{j}e_{i}=0$\textit{\ for }$j>i$%
\textit{. Then, there exists a natural number }$n$\textit{\ such that }$%
e_{i}=0$\textit{\ for all }$i>n$\textit{.\medskip }

\noindent \textbf{Proof. }We denote by $I$ the right ideal $\sum_{i\in
\mathbb{N}}e_{i}S$. By Lemma~3.8, there exists an element $e\in I^{\times
}(S)$ such that $1\equiv _{I}e$ and $I\subseteq \mathrm{Ann}_{r}(e)$. Since $%
1\equiv _{I}e$, we have
\begin{equation}
1+\sum_{k=1}^{m_{1}}e_{i_{k}}s_{i_{k}}=e+\sum_{l=1}^{m_{2}}e_{j_{l}}s_{j_{l}}
\label{id_fin1}
\end{equation}%
for some $s_{i_{k}},s_{j_{l}}\in S$, $k=1,\ldots ,m_{1}$, $l=1,\ldots ,m_{2}$%
. Let $n=\max_{k,l}(i_{k},j_{l})$. Hence, from our assumption we get that $%
e_{t}e_{i_{k}}=0=e_{t}e_{j_{l}}$ for every $t>n$. Whence, multiplying (\ref%
{id_fin1}) by $e_{t}$ on left, we have that $e_{t}=e_{t}e$ and, since $%
e_{t}\in I\subseteq \mathrm{Ann}_{r}(e)$, we have $%
e_{t}=e_{t}e_{t}=e_{t}ee_{t}=0$ for all $t>n$.\textit{\ \ \ \ \ \ }$%
_{\square }\medskip $\textit{\ }

Recall (see, for example, \cite{grillet:cs}) that a commutative monoid $%
(M,+,0)$ is called $\pi $\emph{-}\textit{regular} (or epigroup) if every its
element has a power in some subgroup of $M$. Using Clifford representations
of commutative inverse monoids (see, for example, \cite[Theorem 3.2.1]%
{grillet:cs}), it is easy to show that the last condition is equivalent to
the condition that for any $x\in M$, there exist a natural number $n$ and an
element $y\in M$ such that $nx=nx+y+nx$. A semiring $S$ is called \textit{%
additively }$\pi $\textit{-regular} iff its additive reduction $(S,+,0)$ is
a $\pi $-regular monoid; equivalently, there exist a natural number $n$ and
an element $y\in S$ such that $n1=n1+y+n1$. Also, let $\diamond $ be the
congruence on a semiring $S$ which is defined as follows: $a\diamond b$ $%
\Longleftrightarrow $ $na=b+x$ and $n^{\prime }b=a+x^{\prime }$ for some
natural numbers $n,n^{\prime }\geq 1$ and $x,x^{\prime }\in S$.\medskip

\noindent \textbf{Proposition 3.10} \textit{Let }$S$\textit{\ be a
zerosumfree right CP-semiring,} $S^{\diamond }:=S/\diamond $\textit{, and} $%
\theta _{+}:=$ $\equiv _{I^{+}(S)}$\textit{. Then,}

\textit{(1) }$S$\textit{\ is simultaneously a zeroic and additively} $\pi $%
\textit{-regular semiring;}

\textit{(2)} $S_{S}^{\diamond }\cong $ $I^{+}(S)=eS$ \textit{for some }$e\in
I^{\times }(S)$\textit{;}

\textit{(3) }$\theta _{+}$\textit{\ is the universal congruence.\medskip }

\noindent \textbf{Proof. }By \cite[Lemma 2.2]{il'in:v-s}, $S^{\diamond
}=S/\diamond $ is a nonzero additively idempotent semiring; and, by
Lemma~3.7, there exist an element $e\in I^{\times }(S)$ and an $S$%
-isomorphism $\psi :S_{S}^{\diamond }\longrightarrow eS$ such that $%
1\diamond e$ and $p\psi =1_{S^{\diamond }}$, where $p\colon S\rightarrow
S^{\diamond }$ is the canonical projection. In particular, $1\diamond e$
implies $n1=e+x$ and $n^{\prime }e=1+x^{\prime }$ for some $n,n^{\prime }\in
\mathbb{N}$ and $x,x^{\prime }\in S$. For $S^{\diamond }$ is additively
idempotent, $e\in I^{+}(S)$, we have $e=n^{\prime }e=1+x^{\prime }$ and
\begin{equation*}
n1=e+x=ne+x+ne=n(1+x^{\prime })+x+n(1+x^{\prime })=n1+(nx^{\prime
}+x+nx^{\prime })+n1\text{.}
\end{equation*}%
Thus, $S$ is an additively $\pi $-regular semiring.

Obviously, $eS=\psi (S^{\diamond })\subseteq I^{+}(S)$. Moreover, for each $%
a\in I^{+}(S)$, we get that $p(a)\in S^{\diamond }=p(eS)$, and hence, $%
p(a)=p(b)$ for some $b\in eS$; and hence, $na=b+x$ and $n^{\prime
}b=a+x^{\prime }$ for some natural numbers $n,n^{\prime }\geq 1$ and $%
x,x^{\prime }\in S$. For $a$ and $b$ are additively idempotent elements, $%
a=na=b+x$ and $b=n^{\prime }b=a+x^{\prime }$. Whence $a=b+x=b+b+x=b+a=a+x^{%
\prime }+a=a+x^{\prime }=b$, and hence, $I^{+}(S)=eS$.

By Lemma~3.8, there exists an element $f\in I^{\times }(S)$ such that $1$ $%
\theta _{+}$ $f$ and $I^{+}(S)=eS\subseteq \mathrm{Ann}_{r}(f)$. In fact, $%
f=0$: Indeed, for some $x^{\prime }\in I^{+}(S)$, we have $%
0=fe=f(1+x^{\prime })=f+fx^{\prime }$, and, for $S$ is zerosumfree, $f=0$
and $\theta _{+}$ is the universal congruence. Finally, $1$ $\theta _{+}$ $0$
implies $1+a=b$ for some $a,b\in I^{+}(S)$ and, hence, $1+z=z$ for $z=a+b$
and, therefore, $S$ is a zeroic semiring.\textit{\ \ \ \ \ \ }$_{\square
}\medskip $\textit{\ }

Now, applying Proposition~3.5, Corollary 3.6 and \cite[Theorem 1.2.8]%
{lam:afcinr}, we conclude this section with the following important result
regarding the \textquotedblleft structure\textquotedblright\ of
CP-semirings.\medskip

\noindent \textbf{Theorem 3.11} \textit{A semiring }$S$\textit{\ is a right
(left) CP-semiring iff }$S=R\oplus T$\textit{, where }$R$\textit{\ and }$T$%
\textit{\ are a semisimple ring and a zeroic additively }$\pi $\textit{%
-regular right (left) CP-semiring, respectively.\medskip }

\noindent \textbf{Proof. }$\Longrightarrow $. Let $S$ be a right
CP-semiring, $\equiv _{V(S)}$ the Bourne congruence on $S$ by the ideal $V(S)$%
. It is clear that $\overline{S}:=S/\equiv _{V(S)}$ is a zerosumfree semiring. By
Proposition~3.5, $\overline{S}$ is a right CP-semiring, too. Therefore, by
Proposition~3.10, $\overline{S}$ is a zeroic additively $\pi $-regular
semiring. Thus, applying \cite[Proposition~2.9]{il'in:v-s}, we get $%
S=R\oplus T$, where $R$ is a ring and $T$ is a semiring isomorphic to $%
\overline{S}$. Once again, using Proposition~3.5, one sees that the ring $R$
is a right CP-semiring, and hence, by \cite[Theorem 1.2.8]{lam:afcinr}, $R$
is a semisimple ring.

$\Longleftarrow $. This follows from Corollary 3.6 and \cite[Theorem 1.2.8]%
{lam:afcinr}.\textit{\ \ \ \ \ \ }$_{\square }\medskip $\textit{\ }

\section{Semisimple, Gelfand, subtractive, and anti-bounded, CP-semirings}

In this section, together with some other results, we provide full
descriptions of semisimple, Gelfand, subtractive, and anti-bounded,
semirings the full c-diagram of every semimodule over which is projective,
\textit{i.e.}, the descriptions of semirings of those classes that are also
CP-semirings.

A) Let us start with semisimple (semi)rings. First, applying \cite[Theorem
1.2.8]{lam:afcinr} and \cite[Corollary 21.9]{andfull:racom}, we immediately
observe that for a ring $R$ to be a right (left) CP-ring is a Morita
invariant property. However, taking into consideration \cite[Theorem 5.14]%
{kat:thcos}, we will see in Proposition 4.2 that is not true for semirings
in general. To show that, we will need the following quite useful
observation that matrix semirings over a semiring $S$ and the semiring $S$
itself have isomorphic congruence lattices, namely\medskip

\noindent \textbf{Theorem 4.1} (\textit{cf.} \cite[Theorem 9.1.9]{rs:tqtofs}%
) \textit{Let }$S$\textit{\ be a semiring and }$T=M_{n}(S)$\textit{\ a
matrix semiring over }$S$\textit{.} \textit{For any congruence }$\theta $%
\textit{\ on }$S$\textit{, one defines a congruence }$\Theta $\textit{\ on }$%
T$\textit{\ as follows: for any }$A=(a_{ij}),B=(b_{ij})\in T,$ $A$%
\textit{\ }$\Theta $\textit{\ }$B$\textit{\ }$\Longleftrightarrow $\textit{\
}$a_{ij}$\textit{\ }$\theta $\textit{\ }$b_{ij}$\textit{\ for all }$i,j$%
\textit{. Then, the map }$\chi \colon Cong(S)\ni \theta \mapsto \Theta \in
Cong(T)$\textit{\ is an isomorphism of lattices of congruences, and }$%
T/\Theta \cong M_{n}(S/\theta )$\textit{.\medskip }

\noindent \textbf{Proposition 4.2} \textit{Let a semiring }$S$\textit{\ be
not a ring, then matrix semirings }$M_{n}(S)$\textit{\ are not right (left)
CP-semirings for any }$n\geq 3$\textit{.\medskip }

\noindent \textbf{Proof. }Consider the congruences $\diamond _{S}$ and $%
\diamond _{T}$ on the semirings $S$ and $T=M_{n}(S)$, respectively, defined
as was done before Proposition~3.10. For Theorem 4.1, we readily see that $%
\diamond _{T}=\chi (\diamond _{S})$, and hence, $T/\diamond _{T}\cong
M_{n}(S/\diamond _{S})$. Then, for Proposition~3.5 and~\cite[Lemma~2.2]%
{il'in:v-s} and without loss of generality, one may assume that $S$ is an
additively idempotent semiring.

Let $A\in T$ be the matrix $A=(a_{ij})$ such that $a_{ij}=0$ for $i=j$, and $%
a_{ij}=1$ otherwise. It is sufficient to show that the cyclic right $T$%
-semimodule $AT$ is not projective. Indeed, suppose that it is not the case,
then there exists an injective $T$-homomorphism $\psi \colon AT\rightarrow T$
such that $\alpha \psi =1_{AT}$, where the surjective $T$-homomorphism $%
\alpha :T\longrightarrow AT$ is defined by $\alpha (X)=AX$. Consider the
matrix $B=\psi (AE_{11})\in T$, $B=(b_{ij})$. Since $E_{11}^{2}=E_{11}$, we
get
\begin{equation*}
B=\psi (AE_{11})=\psi (AE_{11}^{2})=\psi (AE_{11})E_{11}=BE_{11}\text{.}
\end{equation*}%
It implies $b_{ij}=0$ for all $j\neq 1$. Moreover, $AE_{11}=\alpha \psi
(AE_{11})=\alpha (B)=AB$, and therefore, $0=(E_{11}A)_{11}=(AB)_{11}=%
\sum_{i=2}^{n}b_{i1}$. From this and the additive idempotentness of $S$, it
follows that $b_{i1}=0$ for all $i>1$.

Finally, $1=(AE_{11})_{21}=(AB)_{21}=a_{21}b_{11}=b_{11}$. So, $\psi
(AE_{11})=B=E_{11}$. Similarly, $\psi (AE_{ii})=E_{ii}$ for all $i=1,\ldots
,n$. The latter implies $\psi (A)=\sum_{i=1}^{n}\psi
(AE_{ii})=\sum_{i=1}^{n}E_{ii}=E$. It is easy to see that $%
A(E+E_{12})=A(E+E_{13})$, hence, $E+E_{12}=\psi (A(E+E_{12}))=\psi
(A(E+E_{13}))=E+E_{13}$, what leads us to a contradiction.\textit{\ \ \ \ \
\ }$_{\square }\medskip $\textit{\ }

As usual, a semiring $S$ is said to be \textit{right (left) semisimple} if
the right (left) regular semimodule is a direct sum of right (left) minimal
ideals. As well known (see, for example, \cite[Theorem 7.8]{hebwei:hoa}, or
\cite[Theorem 4.5]{knt:ossss}), a semiring $S$ is (right, left) semisimple
iff $S\cong M_{n_{1}}(D_{1})\times \ldots \times M_{n_{r}}(D_{r})$, where $%
M_{n_{i}}(D_{i})$ is the semiring of $n_{i}\times n_{i}$-matrices over a
division semiring $D_{i}$ for each $i=1,\ldots ,r$. Using this observation,
in order to describe semisimple CP-semirings, our considerations should be
naturally reduced to the ones of the matrix CP-semirings over division
semirings; and, therefore, from Proposition~4.2, we obtain the
following\medskip\ result.

\noindent \textbf{Proposition 4.3} \textit{A matrix semiring }$S=M_{n}(D)$%
\textit{\ over a division semiring }$D$\textit{\ is a right (left)
CP-semiring iff }$D$\textit{\ is a division ring, or }$D\cong B$\textit{\
and }$n=1,2$\textit{.\medskip }

\noindent \textbf{Proof. }$\Longrightarrow $. Let $S=M_{n}(D)$ be a right
CP-semiring. If $n\geq 3$, by Proposition~4.2, $D$ is a ring. Hence, we need
to consider only the case with $n=1,2$.

Let $D$ be a proper division semiring. Then, the partition $D=\{0\}\cup
D\backslash \{0\}$ defines a congruence $\tau $ on $D$. Obviously, $D/\tau
\cong \mathbf{B}$; hence, by Theorem 4.1, we have $S/\chi (\tau )\cong M_{n}(%
\mathbf{B})$. On the other hand, by Lemma~3.7, the right $S$-semimodule $%
\overline{S}:=S/\chi (\tau )$ is isomorphic to $AS$ via some $S$-isomorphism
$\psi $ such that $\psi (\overline{E})=A$. Clearly, for any non-zero $d\in D$%
, we have $Ed$ $\chi (\tau )$ $E$, hence, $Ad=\psi (\overline{E}d)=\psi (%
\overline{E})=A$. Since $A\neq 0$ and $D$ is a division semiring, the
equality $Ad=A$ implies $D=\{0,1\}\cong \mathbf{B}$.

$\Longleftarrow $. If $D$ is a division ring, then $S$ is a semisimple ring,
and hence, the statement is trivial. Also, it is clear that $\mathbf{B}$ is
a right CP-semiring. Therefore, we need to show only that $M_{2}(\mathbf{B})$
is a right CP-semiring, too. The proof of \cite[Proposition 4.9]%
{aikn:ovsasaowcsai} serves in our case as well; and just for the reader's
convenience, we briefly sketch it here: Namely, we will use the equivalence
of the semimodule categories $\mathcal{M}_{M_{2}(\mathbf{B})}$ and $\mathcal{%
M}_{\mathbf{B}}$ established in \cite[Theorem 5.14]{kat:thcos}: $F:$ $%
\mathcal{M}_{M_{2}(\mathbf{B})}\rightleftarrows $ $\mathcal{M}_{\mathbf{B}%
}:G $, $F(A)=AE_{11}$ and $G(B)=B^{n}$, where $E_{11}\ $is the matrix unit
in $M_{2}(\mathbf{B})$. Let $M$ be a cyclic right $M_{2}(\mathbf{B})$%
-semimodule with a surjective $M_{2}(\mathbf{B})$-homomorphism $f:M_{2}(%
\mathbf{B})\twoheadrightarrow M.$ By \cite[Lemma 4.7]{kn:meahcos}, $F(f):%
\mathbf{B}^{2}\cong F(M_{2}(\mathbf{B}))\twoheadrightarrow F(M)$ is a
surjective $\mathbf{B}$-homomorphism and, hence, there exists the natural
congruence $\equiv _{F(f)}$ on $\mathbf{B}_{\mathbf{B}}^{2}$ such that $%
\mathbf{B}^{2}/\equiv_{F(f)}$ $\cong F(M)$. It is clear that $\{(0,0)\},$ $%
\mathbf{B},$ $\{(0,0),(1,0),(1,1)\}$ and $\mathbf{B}^{2}$ are, up to
isomorphism, the only quotient semimodules of $\mathbf{B}^{2}$ which are
projective by \cite[Theorem 5.3]{hk:tcos}. Whence, $F(M)$ is a projective
right $\mathbf{B}$-semimodule too, and therefore, using \cite[Lemma 4.10]%
{kn:meahcos}, $M\cong G(F(M))$ is a projective right $M_{2}(\mathbf{B})$%
-semimodule as well.\textit{\ \ \ \ \ \ }$_{\square }\medskip $\textit{\ }

Applying Corollary 3.6 and Proposition 4.3, one immediately obtains a full
description of semisimple right (left) CP-semirings, namely\medskip

\noindent \textbf{Theorem 4.4 }\textit{A semisimple semiring }$S$\textit{\
is a right (left) CP-semiring iff }$S\cong S_{1}\times \ldots \times S_{r}$%
\textit{, where every}$\ S_{i},i=1,\ldots ,r$\textit{, is either an Artinian
simple ring, or isomorphic to }$M_{n}(\mathbf{B})$ \textit{with }$n=1,2$%
\textit{.\medskip }

B) We say that a semiring $S$ is \textit{right (left) Gelfand} (see\textit{\
}also \cite[page 56]{golan:sata}) if for every element $s\in S$, the element
$1+s$ has a right (left) inverse. In this subsection, we present a full
description of Gelfand semirings that simultaneously are CP-semirings. In
order to provide this description, we need first to justify the following
useful observations.\medskip

\noindent \textbf{Lemma 4.5}\textit{\ Let }$I\subseteq S$\textit{\ be a
right ideal of a right Gelfand semiring }$S$\textit{. Then, }$\equiv _{I}$%
\textit{\ is the universal congruence iff }$I=S$\textit{.\medskip }

\noindent \textbf{Proof. }$\Longrightarrow $. Let $\equiv _{I}$ be the
universal congruence. Particularly, we have $1\equiv _{I}0$, that is, $1+a=b$
for some elements $a,b\in I$. Since $S$ is a right Gelfand semiring, $1+a$
has a right inverse $c\in S$, therefore, $1=(1+a)c=bc\in I$, whence $I=S$.

$\Longleftarrow $. This is obvious.\textit{\ \ \ \ \ \ }$_{\square }\medskip
$\textit{\ }

\noindent \textbf{Lemma 4.6}\textit{\ Every right CP- and right Gelfand
semiring is an additively idempotent semiring.\medskip }

\noindent \textbf{Proof. }For a right CP- and right Gelfand semiring $S$, by
Theorem~11, $S=R\oplus T$ where $R$ and $T$ are a semisimple ring and a
zeroic semiring, respectively. Denoting by $1_{R}$ and $1_{T}$ the
multiplicative identities of $R$ and $S$, respectively, we have $%
1_{T}=1-1_{R}$, and hence, $1_{T}$ has a right inverse. The latter obviously
implies that $T=S$, that is, $S$ is a zeroic semiring. Then, by
Proposition~3.10 and Lemma 4.5, $S=I^{+}(S)$.\textit{\ \ \ \ \ \ }$_{\square
}\medskip $\textit{\ }

\noindent \textbf{Lemma 4.7}\textit{\ If }$S$\textit{\ is a right CP- and
right Gelfand semiring, then, for any }$e,f\in I^{\times }(S)$\textit{, the
following statements are true:}

\textit{(1)~}$1+e=1$\textit{;}

\textit{(2)~}$e+f\in I^{\times }(S)$\textit{;}

\textit{(3)~}$eS+fS=(e+f)S$\textit{;}

\textit{(4)~}$eS=fS$\textit{\ }$\Rightarrow $\textit{\ }$e=f$\textit{%
.\medskip }

\noindent \textbf{Proof. }(1) As by Lemma~4.6 $S$ is additively idempotent,
it is easy to see that $1+e\in I^{\times }(S)$. The latter, since $1+e$ has
a right inverse, right away implies that $1+e=1$.

(2) For (1), $(e+f)^{2}=e+ef+fe+f=e(1+f)+f(1+e)=e+f.$

(3) Obviously, $(e+f)S\subseteq eS+fS$. On the other hand, $%
(e+f)e=e+fe=(1+f)e=e$ and similarly $(e+f)f=f$. Hence, $%
eS+fS=(e+f)eS+(e+f)fS\subseteq (e+f)S$.

(4) Assume that $eS=fS$. Then, in particular, we get $e=fe$, and hence, $%
f+e=f+fe=f(1+e)=f$; similarly, $e+f=e$, and therefore, $e=f$.\textit{\ \ \ \
\ \ }$_{\square }\medskip $\textit{\ }

\noindent \textbf{Proposition 4.8}\textit{\ Let }$I$\textit{\ be a right
ideal of a right CP- and right Gelfand semiring }$S$\textit{. Then, there
exist elements }$e,f\in I^{\times }(S)$\textit{\ such that }$e+f=1$\textit{,
}$ef=fe=0$\textit{, and }$I=fS$\textit{. In particular, }$S$\textit{\ is a
right subtractive semiring, and }$I$\textit{\ is a direct summand of }$S_{S}$%
\textit{.\medskip }

\noindent \textbf{Proof. }For $I=S$, or $I=0$, the statement is trivial,
we assume that $I\subset S$ and $I\neq 0$. Then, by Lemma~4.5, the
congruence $\equiv _{I}$ is not universal, and hence, $\overline{S}%
:=S/\equiv _{I}$ is a non-zero semimodule. Therefore, by Lemma~3.8, there
exists an element $e\in I^{\times }(S)$ such that $1\equiv _{I}e$ and $%
I\subseteq \mathrm{Ann}_{r}(e)$. In particular, we have $1+a=e+b$ for some $%
a,b\in I$. Since $1+a$ has a right inverse $c\in S$, we obtain $%
1=(1+a)c=(e+b)c=ec+bc$; and hence, $e=e(bc+ec)=e^{2}c+ebc=ec$, since $e\in
I^{\times }(S)$ and $bc\in I\subseteq \mathrm{Ann}_{r}(e)$. So, $1=ec+bc=e+bc
$, and multiplying it by $bc$ on right, one obtains $bc=ebc+(bc)^{2}=(bc)^{2}
$, \emph{i.e.}, $bc\in I^{\times }(S)$.

Denoting $bc\in I$ by $f$, we have $e+f=1$ and $ef=0$. Also, $e+f=1$ implies
$x=ex+fx=fx$ for all $x\in I$, so, $I=fS$. Clearly, $f\not\in \{0,1\}$, and
hence, $J=eS$ is a right ideal such that $J\subset S$ and $J\neq 0$.
Applying the above reasoning to $J$ and keeping in mind Lemma~4.7 (4), one
gets $g+e=1$, $ge=0$ for some $g\in I^{\times }(S)$. Therefore, $%
f=(g+e)f=gf+ef=gf=g(f+e)=g$ and $fe=ge=0$. Thus, the idempotents $e$ and $f$
are mutually orthogonal and, applying the Pierce's decomposition to the
regular semimodule $S_{S}$, we have $S_{S}=fS\oplus eS=I\oplus J$.

Finally, for $I$ is a direct summand of $S_{S}$, it is subtractive. So, the
rest is obvious.\textit{\ \ \ \ \ \ }$_{\square }\medskip $\textit{\ }

The following theorem, providing a full description of right Gelfand
semirings that are right CP-semirings as well, also solves Problem 2 left
open in \cite{aikn:ovsasaowcsai}.\medskip

\noindent \textbf{Theorem 4.9 }\textit{A right (left) Gelfand semiring }$S$%
\textit{\ is a right (left) CP-semiring iff }$S$\textit{\ is a finite
Boolean algebra.} \medskip

\noindent \textbf{Proof. }$\Longrightarrow $. Assume that $S$ is
simultaneously a right CP- and a right Gelfand semiring. By Proposition 4.8
and \cite[Theorem 4.4]{knt:ossss}, $S$ is a semisimple right CP-semiring.
Then, using Theorem~4.4 and Lemma~4.6, one immediately gets that $S\cong
S_{1}\times \ldots \times S_{r}$, where each $S_{i}$, $i=1,\ldots ,r$, is
isomorphic to $M_{n}(\mathbf{B})$ and $n=1,2$. However, as $S$ is a right
Gelfand semiring, each $S_{i}$, $i=1,\ldots ,r$, is isomorphic to $\mathbf{B}
$, that is, $S$ is a finite Boolean algebra.

$\Longleftarrow $. Let $S$ be a finite Boolean algebra. Then, $S$ is a
direct sum of finitely many copies of $\mathbf{B}$. Using Theorem~4.4, one
ends the proof.\textit{\ \ \ \ \ \ }$_{\square }\medskip $

C) In the following result, using Theorems 3.11 and 4.9, we obtain a full
description of right subtractive right CP-semiring.\medskip

\noindent \textbf{Theorem 4.10 }\textit{A right (left) subtractive semiring }%
$S$\textit{\ is a right (left) CP-semiring iff }$S=R\oplus T$\textit{, where
}$R$\textit{\ and }$T$\textit{\ are a semisimple ring and a finite Boolean
algebra, respectively. \medskip }

\noindent \textbf{Proof. }$\Longrightarrow $. Let $S$ be a right
subtractive, right CP-semiring. By Theorem~3.11, $S=R\oplus T$, where $R$
and $T$ are a semisimple ring and a zeroic right CP-semiring, respectively.
For \cite[Lemma 4.7]{knt:ossss}, $T$ is a right subtractive semiring, too.
Since by Proposition~3.10~(3), the congruence $\theta _{+}$ on $T$ is
universal, we have $1$ $\theta _{+}$ $0$, \emph{i.e.}, $1+a=b$ for some $%
a,b\in I^{+}(T)$. Whence, $1\in I^{+}(T)$ since $T$ is right subtractive,
and hence, $T=I^{+}(T)$. Then, for each $x\in T$, we get $1+(1+x)=1+x$; so,
for $T$ is right subtractive, $1\in (1+x)T$. Therefore, $1+x$ has a right
inverse in $T$, that is, $T$ is a right Gelfand semiring. So, applying
Theorem~4.9, one sees right away that $T$ is a finite Boolean algebra.

$\Longleftarrow $. Assume that $S=R\oplus T$, where $R$ and $T$ are a
semisimple ring and a finite Boolean algebra, respectively. Obviously, $R$
and $T$ are right subtractive semirings; hence, by \cite[Lemma 4.7]%
{knt:ossss}, $S$ is right subtractive, too. Next, using Theorems~3.11
and~4.9, one ends the proof.\textit{\ \ \ \ \ \ }$_{\square }\medskip $

As was mentioned earlier, the concepts of `congruence-simpleness' and
`ideal-simpleness' for (even for subtractive) semirings are not the same
(see, for example, \cite[Examples 3.8]{knz:ososacs}, or \cite[Theorems 3.7
and 4.5]{knt:mosssparp}). Also, as were shown in \cite[Theorem 14.1]%
{bashkepka:css} and \cite[Corollary 4.4]{knt:mosssparp}, respectively, for
finite commutative and left (right) subtractive semirings, these concepts
coincide. As a corollary of Theorem 4.10, our next result, solving \cite[%
Problem 3]{knt:mosssparp} and \cite[Problem 3]{knt:mosssparp}, shows that
for subtractive semirings these two concepts coincide as well.\medskip

\noindent \textbf{Corollary 4.11 }\textit{For a right subtractive semiring }$%
S,$\textit{\ the following conditions are equivalent:}

\noindent\ \ \ \ \ \textit{(1) }$S$\textit{\ is a congruence-simple right
CP-semiring;}

\textit{(2) }$S$\textit{\ is an ideal-simple right CP-semiring;}

\textit{(3) }$S\cong M_{n}(D)$\textit{\ for some division ring }$D,$\textit{%
\ or }$S\cong \mathbf{B}$\textit{.\medskip }

\noindent \textbf{Proof. }(1) $\Longrightarrow $ (2). This follows
immediately from \cite[Proposition 4.4]{knt:mosssparp}.

(2) $\Longrightarrow $ (3). Assume that $S$ is an ideal-simple right
CI-semiring. By Theorem 4.10, $S$ is either a semisimple ring, or a finite
Boolean algebra. If $S$ is a semisimple ring, then, since $S$ is
ideal-simple, $S\cong M_{n}(D)$ for some division ring $D$ and $n\geq 1$.
Otherwise, $S$ is a finite Boolean algebra. By \cite[Theorem 3.7]%
{knt:mosssparp}, $S$ is a proper division semiring, and, since $S$ is a
finite Boolean algebra, $S$ is just the Boolean semifield $\mathbf{B}$.

(3) $\Longrightarrow $ (1). This follows immediately from Theorem~3.11,
Theorem 4.10 and \cite[Theorem 4.5]{knt:mosssparp}.\textit{\ \ \ \ \ \ }$%
_{\square }\medskip $

D) Following \cite{aikn:ovsasaowcsai}, a semiring $S$ is said to be \textit{%
anti-bounded} if $S=V(S)\cup \{1+s\,|\,s\in S\}$. And we conclude the
current section with a full description of anti-bounded CP-semirings ---
what constitutes one of the main results of these section and paper. But
first let us consider some important examples of anti-bounded
CP-semirings.\medskip

\noindent \textbf{Facts 4.12 }

(1) \textit{The semiring} $\mathbf{B}_{3}$\textit{, defined on the chain }$%
0<1<2$\textit{, with the addition }$x+y:=x\vee y$\textit{\ and multiplication%
}%
\begin{equation*}
xy:=\left\{
\begin{array}{lcl}
0\text{\textit{,}} &  & \text{\textit{if} }x=0\,\,\text{\textit{or}}\,\,y=0
\\
&  &  \\
x\vee y\text{\textit{,}} &  & \text{\textit{otherwise}}%
\end{array}%
\right. \text{\textit{,}}
\end{equation*}%
\textit{is a CP-semiring.}

\textit{(2) The semiring }$B(3,1)=(\{0,1,2\},\oplus ,\odot )$\textit{\ with
the operations }$a\oplus b\overset{def}{=}\min (2,a+b)$\textit{\ and }$%
a\odot b\overset{def}{=}\min (2,ab)$\textit{\ is a CP-semiring.\medskip }

\noindent \textbf{Proof. }(1) As was shown in \cite[Fact 4.11]%
{aikn:ovsasaowcsai}, up to isomorphism, there are only two nonzero cyclic $%
\mathbf{B}_{3}$-semimodules, namely $\{0,2\}$ and $\mathbf{B}_{3}$; and it
is easy to see that these semimodules are retracts of the regular semimodule
$\mathbf{B}_{3}$. So, they are projective, and $\mathbf{B}_{3}$ is a
CP-semiring.

(2) This fact can be verified in a similar to (1) fashion.\textit{\ \ \ \ \
\ }$_{\square }\medskip $

For the reader's convenience, we remind here another important class of
anti-bounded semirings quite naturally arising from rings and originally
introduced in \cite[Example 4.16]{aikn:ovsasaowcsai}. Let $R=(R,{+},{\cdot }%
,e,1)$ be an arbitrary ring with zero $e$ and unit $1$. Let $T:=R\cup \{0\}$
and extend the operations on $R$ to $T$ by setting $0+t=t=t+0$ and $0\cdot
t=0=t\cdot 0$ for all $t\in T$. Clearly, $(T,+,\cdot ,0,1)$ is a zerosumfree
semiring. Now, extend the semiring structure on $T$ to a semiring structure
on $\mathrm{Ext}(R):=T\cup \{\infty \}=R\cup \{0,\infty \}$, where $\infty
\notin T$, by setting $x+\infty =\infty =\infty +\infty =\infty +x$ and $%
x\cdot \infty =\infty =\infty \cdot \infty =\infty \cdot x$ for all $x\in R$%
, and $0\cdot \infty =0=\infty \cdot 0$. It is easy to see that $(\mathrm{Ext%
}(R),+,\cdot ,0,1)$ is, indeed, an anti-bounded zerosumfree semiring. In a
similar fashion, one can naturally extend the structure of every right $R$%
-module $M$ to a structure of an $\mathrm{Ext}(R)$-semimodule on $\mathrm{Ext%
}(M)$.\medskip

\noindent \textbf{Proposition 4.13}\textit{\ For a ring }$R$\textit{, the
semiring }$\mathrm{Ext}(R)$ \textit{is a right (left) CP-semiring iff }$R$%
\textit{\ is a semisimple ring.\medskip }

\noindent \textbf{Proof. }By \cite[Proposition 4.17]{aikn:ovsasaowcsai}, up
to isomorphism, $\{0\}$, $\{0,\infty \}$, and $\{\mathrm{Ext}(\overline{R})$
$|$ $\overline{R}=R/I\}$, where $I$ is a right ideal of $R$, are all cyclic
right $\mathrm{Ext}(R)$-semimodules. Obviously, $\{0\}$ and $\{0,\infty \}$
are retracts of $Ext(R)_{Ext(R)}$, and hence, they are projective. Also, one
may readily verify that a cyclic right $R$-module $M$ is a retract of $R_{R}$
iff the cyclic right $Ext(R)$-semimodule $Ext(M)$ is a retract of $%
Ext(R)_{Ext(R)}$. Using these observations and Theorem 3.3, one ends the
proof.\textit{\ \ \ \ \ \ }$_{\square }\medskip $

\noindent \textbf{Proposition 4.14 }\textit{A nonzero additively idempotent
anti-bounded semiring }$S$\textit{\ is a right (left) CP-semiring iff} $%
S\cong \mathbf{B}$\textit{, or }$S\cong \mathbf{B}_{3}$\textit{.\medskip }

\noindent \textbf{Proof. }$\Longrightarrow $. First note that $(S,{+})$ is
an upper semilattice such that $0<1\leq s$ for any nonzero $s\in S$.
Obviously, $S$ is an entire (\emph{i.e.}, has no zero divisors) semiring;
and hence, the partition $S=\{0\}\cup S\backslash \{0\}$ defines a
congruence $\tau $ on $S_{S}$. By Lemma~3.7, there exists an $S$-isomorphism
$\psi \colon \overline{S}:=S/\tau \rightarrow zS$ for some $z\in I^{\times
}(S)$ and $\psi (\overline{1})=z$. Particularly, for every nonzero $a\in S$,
we have $\overline{1}a=\overline{1}$, and hence, $za=z$. Clearly, as $z\neq
0 $, for some $x\in S$, we get $z=1+x\geq 1$ and, therefore, $z=za\geq a$
and $z+a=z$. Thus, $z$ is an infinite element of $S$. Obviously, if $z=1$,
then $S\cong \mathbf{B}$.

Now consider the remaining case when $z>1$. Fix $a\in S$ and $a>1$. Clearly,
$a+s\geq a>1$ for all $s\in S$. Moreover, if $s\neq 0$, then $s=1+x$ for
some $x\in S$; and hence, $as=a(1+x)=a+ax\geq a>1$. Therefore, the partition
$S=\{0\}\cup \{1\}\cup S\backslash \{0,1\}$ defines a congruence $\sigma $
on $S_{S}$. Again, by Lemma~3.7, there exists an injective $S$-homomorphism $%
\varphi \colon \tilde{S}:=S/\sigma \rightarrow S$ such that $\alpha \varphi
=1_{\tilde{S}}$, where $\alpha \colon S\rightarrow \tilde{S}$ is the
canonical projection. In particular, from the latter it immediately follows
that $\varphi (\tilde{1})=1$. Then, $\varphi (\tilde{z})=\varphi (\tilde{1}%
)z=z$. Furthermore, for any $a\in S$, $a>1$, we have $\tilde{a}=\tilde{z}$,
hence, $z=\varphi (\tilde{z})=\varphi (\tilde{a})=\varphi (\tilde{1})a=a$.
Thus, $S$ has actually only three elements: $0$, $1$, and $z$; that is, $%
S\cong \mathbf{B}_{3}$.

$\Longleftarrow $. This immediately follows from Facts 4.12 (1) and
Theorem~4.4.\textit{\ \ \ \ \ \ }$_{\square }\medskip $

\noindent \textbf{Proposition 4.15 }\textit{A zerosumfree anti-bounded
semiring }$S$\textit{\ is a right (left) CP-semiring iff }$S\cong \mathbf{B}$%
\textit{, or }$S\cong \mathbf{B}_{3}$\textit{, or }$S\cong B(3,1)$\textit{,
or }$S\cong $ $\mathrm{Ext}(R)$ \textit{for some nonzero semisimple ring }$R$%
\textit{.\medskip }

\noindent \textbf{Proof. }$\Longrightarrow $. By Proposition~3.5, the
quotient semiring $S^{\diamond }=S/\diamond $ is a right CP-semiring as
well. Then, since $S^{\diamond }$ is an additively idempotent semiring, by
Proposition~4.14, either $S^{\diamond }\cong \mathbf{B}$, or $S\cong \mathbf{%
B}_{3}$. Furthermore, by Proposition~3.10~(2), $S_{S}^{\diamond }\cong
I^{+}(S)$ and, hence, $I^{+}(S)=\{0,\infty \}$, or $I^{+}(S)=\{0,e,\infty \}$%
, and $I^{+}(S)$ possesses an infinite element. By Proposition~3.10~(3), $%
\theta _{+}$ is the universal congruence on $S$; and therefore, for every $%
a\in S$, we have $a$ $\theta _{+}$ $0$ and there exist elements $x,y\in
I^{+}(S)$ such that $a+x=y$. Whence, $a+\infty =a+x+\infty =y+\infty =\infty
$. Thus, $\infty $ is the infinite element of $S$.

Now suppose that $S$ is not additively regular, that is, $1+x+1\neq 1$ for
all $x\in S$. Then, the partition $S=\{0\}\cup \{1\}\cup S\backslash \{0,1\}$
defines a congruence $\sigma $ on $S_{S}$. Indeed, if $a\not\in \{0,1\}$,
then $a=1+c$ for some nonzero $c\in S$. It implies $a=1+x+1$ for some $x\in
S $. Therefore, for any $s\neq 0$, we get $s=1+y$ for suitable $y\in S$;
and, hence, $a+s=1+x+1+s\not\in \{0,1\}$ and $as=a(1+y)=a+ay=1+x+1+ay\not\in
\{0,1\}$. By Lemma~3.7, there exists a $S$-monomorphism $\psi \colon
\overline{S}:=S/\sigma \rightarrow S$ such that $\alpha \psi =1_{\overline{S}%
}$, where $\alpha \colon S\rightarrow \overline{S}$ is the canonical
projection. In particular, as was shown in the proof of Proposition~4.14,
the latter implies $\psi (\overline{1})=1$, $\psi (\overline{\infty })=\psi (%
\overline{1})\infty =\infty $; so, for any $a\in S\backslash \{0,1\}$, we
have $\overline{a}=\overline{\infty }$ and $\infty =\psi (\overline{\infty }%
)=\psi (\overline{a})=\psi (\overline{1})a=a$. Thus, $S$ has, in fact, only
three elements: $0$, $1$, and $\infty $; and hence, $1+1=\infty $. Thus, $%
S\cong B(3,1)$.

Now suppose $S$ is additively regular. Then, using Clifford representations
of commutative inverse monoids (see, for example, \cite[Theorem 3.2.1]%
{grillet:cs}), the additive reduct $(S,+,0)$ is a disjoint union of abelian
groups $(G_{x},+,x)$, $x\in I^{+}(S)$. Clearly, we have $G_{0}=\{0\}$ and $%
G_{\infty }=\{\infty \}$. Hence, if $I^{+}(S)=\{0,\infty \}$, then we get $%
S=I^{+}(S)\cong \mathbf{B}$. Finally, if $I^{+}(S)=\{0,e,\infty \}$, then
one may easily see that the subset $R:=G_{e}\subset S$ is closed under
addition and multiplication, therefore, $(R,+,\cdot ,e,1)$ is a ring. If $R$
is zero, then $S=I^{+}(S)\cong \mathbf{B}_{3}$. Otherwise, one has $S=Ext(R)$
and, therefore, $R$ is a semisimple ring by Proposition~4.13.

$\Longleftarrow $. It follows from Propositions~4.13,~4.14, and Facts 4.12.%
\textit{\ \ \ \ \ \ }$_{\square }\medskip\medskip$

Applying Theorem~3.11 and Proposition~4.15, we have a full description of
anti-bounded CP-semirings, namely:\medskip

\noindent \textbf{Theorem 4.16 }\textit{An anti-bounded semiring }$S$\textit{%
\ is a right (left) CP-semiring iff }$S$\textit{\ is one of the following
semirings:}

\textit{(1) }$S$\textit{\ is a semisimple ring;}

\textit{(2)} $S\cong \mathbf{B}$\textit{, or }$S\cong \mathbf{B}_{3}$\textit{%
, or }$S\cong B(3,1)$\textit{, or }$S\cong $ $\mathrm{Ext}(R)$ \textit{for
some nonzero semisimple ring }$R$\textit{;}

\textit{(3) }$S=R\oplus T$\textit{, where }$R$\textit{\ is a semisimple ring
and }$T$\textit{\ is isomorphic to} $\mathbf{B}$\textit{, or }$\mathbf{B}%
_{3} $\textit{, or }$B(3,1)$\textit{, or }$S\cong $ $\mathrm{Ext}(R')$
\textit{for some nonzero semisimple ring }$R'$.\medskip

Using Theorem~4.16, it is easy to see that the concepts of
congruence-simpleness and ideal-simpleness for anti-bounded right (left)
CP-semirings are the same --- these semirings are isomorphic either to
matrix rings over division rings, or to the Boolean semifield $\mathbf{B}$.
In this connection, however, we conclude this section by presenting a more
general result characterizing anti-bounded semirings for which these two
concepts of simpleness coincide and, therefore, by solving \cite[Problem 3]%
{knt:mosssparp} and \cite[Problem 5]{knz:ososacs} in the class of
anti-bounded semirings.\medskip

\noindent \textbf{Theorem 4.17 }\textit{For an anti-bounded semiring }$S$%
\textit{\ the following conditions are equivalent:}

\textit{(1) }$S$\textit{\ is congruence-simple;}

\textit{(2) }$S$\textit{\ is ideal-simple;}

\textit{(3) }$S$\textit{\ is a simple ring, or} $S\cong \mathbf{B}$\textit{.}%
\medskip

\noindent \textbf{Proof. }(1) $\Longrightarrow $ (3). Assume that $S$ is a
congruence-simple semiring. By \cite[Proposition 3.1]{zumbr:cofcsswz}, $S$
is either a ring, or an additively idempotent semiring. If $S$ is a ring,
then since $S$ is congruence-simple we see right away that $S$ is a simple
ring.

Now consider the case when $S$ is an additively idempotent semiring, and
hence, its additive reduct $(S,+,0)$ is an upper semilattice with the
partial order relation defined for any $x,y\in S$ as $x\leq y$ $%
\Longleftrightarrow $ $x+y=y$. For $S$ is anti-bounded, $0<1\leq x$ for all $%
0\neq x\in S$, and hence, $S$ is entire. Then, there exist the surjective
semiring homomorphism $\pi :S\longrightarrow \mathbf{B}$ such that $\pi
(0):=0$, and $\pi (x):=1$ for all nonzero $x\in S$ and the corresponding
natural congruence $\equiv _{\pi }$ on $S$. Since $(1,0)\notin $ $\equiv
_{\pi }$ and $S$ is a congruence-simple semiring, $\equiv _{\pi }$ is the
diagonal congruence, and hence, $S\cong \mathbf{B}$.

(2) $\Longrightarrow $ (3). Let $S$ be an ideal-simple semiring. For $V(S)$
is an ideal of $S$, we have $V(S)=S$, or $V(S)=0$. If $V(S)=S$, then $S$ is
a simple ring.

Now let $S$ be a zerosumfree semiring. Then, one sees right away that $%
(1+1)S:=\{s+s\,|\,s\in S\}$ is a nonzero ideal of $S$; and as $S$ is
ideal-simple, $(1+1)S=S$, and hence, $s+s=1$ for some nonzero $s\in S$.
Since $S$ is an anti-bounded semiring, there exists an element $x\in S$ such
that $1=1+x+1+x$, and hence, $0\neq 1+x+x\in I^{+}(S)$. By $S$ is
ideal-simple, $I^{+}(S)=S$, \textit{i.e.}, $(S,+,0)$ is an upper semilattice
with the natural order defined as follows: $\forall $ $x,y\in S$ $(x\leq
y\Longleftrightarrow x+y=y)$. Then, since $I:=\{0\}\cup \{s\in S\,|\,s>1\}$
is obviously an ideal of the ideal-simple semiring $S$, we have that $I=0$,
and therefore, $S=\{0,1\}\cong \mathbf{B}$.

The implications (3) $\Longrightarrow $ (1) and (3) $\Longrightarrow $ (2)
are trivial.\textit{\ \ \ \ \ \ }$_{\square }\medskip $

\section{CP-semirings of endomorphisms of semimodules over Boolean algebras}

In the previous section, there have been obtained the full descriptions of
CP-semirings within widely known and important classes of semirings. In
contrast to that, there exists another quite general approach to the problem
of describing CP-semirings. As is clear, any semiring $S$ can be considered
as a semiring of all endomorphisms $End(M_{T})$ of some semimodule $M_{T}\in
$\textit{\ }$|\mathcal{M}_{T}|$ over a semiring $T$ (for instance, $S$ $%
\cong End(S_{S})$, where $S_{S}\in $\textit{\ }$|\mathcal{M}_{S}|$ is a
regular right semimodule, of course). Therefore, the following general
problem/program sounds quite natural and seems to be very interesting and
uneasy as well: Given a semiring $T$, characterize all semimodules $M_{T}\in
$\textit{\ }$|\mathcal{M}_{T}|$ such that their endomorphism semirings $%
End(M_{T})$ are CP-semirings. In this section, we initiate this program
considering the semiring $T$ to be a Boolean algebra --- namely, for an
arbitrary Boolean algebra $B$, we describe all semimodules $M_{B}\in $%
\textit{\ }$|\mathcal{M}_{B}|$ whose endomorphism semirings $End(M_{B})$ are
CP-semirings.\medskip\ We start with a general observation:

\noindent \textbf{Proposition 5.1 }\textit{Let }$e\in I^{\times }(S)$\textit{%
\ be a nonzero idempotent of a right CP-semiring }$S$\textit{. Then, }$eSe$%
\textit{\ is a right CP-semiring, too.\medskip }

\noindent \textbf{Proof. }Let $S$ be a right CP-semiring and $0\neq e\in
I^{\times }(S)$. Then, there exist the restriction and induction functors $%
Res:$ $\mathcal{M}_{S}\longrightarrow $ $\mathcal{M}_{eSe}$ and $Ind:$ $%
\mathcal{M}_{eSe}\longrightarrow $ $\mathcal{M}_{S}$, respectively, given by
\begin{equation*}
Res(M)=Me,\,\,\,\,Ind(M)=M\otimes _{eSe}eS.
\end{equation*}%
Moreover, it is easy to see that $Ind$ is a left adjoint of $Res$, and there
is a natural bijection $Res\circ Ind\cong Id_{\mathcal{M}_{eSe}}$; and,
hence, $Ind$ preserves colimits and $Res$ preserves limits (see, \textit{e.g.%
}, \cite[Theorem 5.5.1]{macl:cwm}), also one may see that $Res$ preserves
finite colimits as well.

Now, let $M$ be a cyclic right $eSe$-semimodule and $f:eSe\longrightarrow M$
a surjective $eSe$-homomorphism. It is clear that $f$ is a cokernel of its
own kernel pair, and therefore, the $S$-homomorphism $Ind(f):eS\cong
Ind(eSe)\longrightarrow Ind(M)$ is also a cokernel of its own kernel pair
and, hence, surjective, too. For $S$ is a right CP-semiring, the latter
implies that the cyclic $S$-semimodule $Ind(M)$ is projective. Then, since
the functor $Res$ preserves finite colimits, it preserves, in particular,
the projectiveness of semimodules, and, using the isomorphism $M\cong
Res(Ind(M))$, one concludes that the cyclic right $eSe$-semimodule $M$ is
projective and $eSe$ is a right CP-semiring.\textit{\ \ \ \ \ \ }$_{\square
}\medskip $

From now on, unless otherwise stated, let $M=(M,+,0)$ be a right $\mathbf{B}$%
-semimodule and $S:=End(M)$. Obviously, $(M,+,0)$ is an upper semilattice $%
(M,{\vee },{0})$ with the least element $0$ and the natural order defined
as: $a\leq b$ $\Leftrightarrow $ $a+b=b$. For each $m\in M$, let $%
m^{\triangledown }\!:=\!\{x\in M:x\leq m\}$, $m^{\vartriangle }\!:=\!\{x\in
M:m\leq x\}$. Also, for any $a,b\in M$, let $e_{a,b}$ be the map from $M$ to
$M$, defined for any $x\in M$ as follows:
\begin{equation*}
e_{a,b}(x):=\{_{b\text{ otherwise \ }}^{0\text{ if }x\text{ }\leq \text{ }a,}%
\text{.}
\end{equation*}%
It is obvious (also, see \cite[Lemma 2.2]{zumbr:cofcsswz}) that $e_{a,b}\in
S $; moreover, $e_{a,b}\in I^{\times }(S)$ provided $b\not\in
a^{\triangledown }$.\medskip

\noindent \textbf{Lemma 5.2 }\textit{If }$S$\textit{\ is a right
CP-semiring, then }$M$\textit{\ satisfies the ascending chain
condition.\medskip }

\noindent \textbf{Proof. }Let $S$ be a right CP-semiring, and $%
x_{1}<x_{2}<\dots $ an infinite ascending chain in $M$. Then, the set $%
\{f_{i}\in S,\,i\geq 1\}$, where $f_{i}=e_{x_{i},x_{i+1}}$, is a subset of $%
I^{\times }(S)$. Obviously, $f_{i}\neq 0$ for each $i$, and $f_{j}f_{i}=0$
for $j>i$. However, the latter contradicts to Lemma~3.9.\textit{\ \ \ \ \ \ }%
$_{\square }\medskip $

\noindent \textbf{Lemma 5.3 }\textit{If }$S$\textit{\ is a right
CP-semiring, then }$M$\textit{\ is a bounded lattice.\medskip }

\noindent \textbf{Proof. }Obviously, every upper semilattice satisfying the
ascending chain condition has the greatest element, hence, by Lemma~5.2, $M$
is a bounded semilattice. Also, it is easy to see that, for every $a,b\in M$%
, the set $a^{\triangledown }\cap b^{\triangledown }$ is a subsemilattice of
$M$; so, it possesses the greatest element $c$ as well. Clearly, $c$ is the
greatest lower bound of $a$ and $b$, \emph{i.e.}, $c=a\wedge b$. Thus, $M$
is a bounded lattice.\textit{\ \ \ \ \ \ }$_{\square }\medskip $

The following facts, establishing that endomorphism semirings of a \textit{%
pentagon }$\mathit{N}_{5}$\textit{\ or diamond }$\mathit{M}_{3}$ (see, for
instance, \cite[p.79]{gratzer:glt}) as well as of a bounded infinite
descending chain are not CP-semirings, will prove to be useful.\medskip\

\noindent \textbf{Fact 5.4 }\textit{Let }$M$\textit{\ be a lattice
isomorphic to }$M_{3}$\textit{, i.e., }$M=\{0,a,b,c,1\}$\textit{, where }$0$%
\textit{\ and }$1$\textit{\ are its least and greatest elements,
respectively, and }$a$\textit{, }$b$\textit{, }$c$\textit{\ are mutually
incomparable elements. Then, }$End(M)$\textit{\ is not a right
CP-semiring.\medskip }

\noindent \textbf{Proof. }Let $S:=End(M)$ and consider the relation $\theta
_{0}$ on $S$ defined for all $s,s^{\prime }\in S$ as follows: $s$ $\theta
_{0}$ $s^{\prime }$ $\Longleftrightarrow $ $s=ru+v$, $s^{\prime }=ru^{\prime
}+v$ for some $u,v\in S$ and $r,r^{\prime }\in \{e_{0,c},\,e_{0,1}\}$. Then,
let $\theta $ be the congruence on $S_{S}$ generated by the relation $\theta
_{0}$, \textit{i.e.}, $\theta $ is the transitive closure of $\theta _{0}$.

Now notice that if $u\neq 0$ and $r\in \{e_{0,c},\,e_{0,1}\}$, then there
exist at least two different elements $x,y\in \{a,b,c\}$ such that $u(x)\neq
0$ and $u(y)\neq 0$: Indeed, if, for instance, $u(a)=u(b)=0$, then $%
u(1)=u(a+b)=u(a)+u(b)=0$ and, consequently, $u=0$. Thus, we have $u(x)\neq 0$
and $u(y)\neq 0$ and, hence, $ru(x)=ru(y)\neq 0$. From the latter it is easy
to see that the unit $1_{S}$ cannot be represented in the form $ru+v$ with $%
r\in \{e_{0,c},e_{0,1}\}$ and $u\neq 0$. Whence, $1_{S}$ $\theta _{0}$ $s$
implies $1_{S}=s$, and therefore, $1_{S}$ $\theta $ $s$ implies $1_{S}=s$.

We claim that $\overline{S}:=S/\theta $ is not a retract of $S_{S}$. Assume
that is not a case. Then, by Lemma~3.7, there exists an element $e\in
I^{\times }(S)$ such that $\overline{S}\cong eS$ and $1_{S}$ $\theta $ $e$.
Whence, $1_{S}=e$, and hence, $\overline{S}\cong eS=S$ and $|\overline{S}%
|=|S|$. However, since $M$ is finite, so is $S$, and $|\overline{S}|<|S|$
since $\theta $ $\neq $ $\vartriangle _{S}$. Thus, $\overline{S}$ is not a
retract of $S_{S}$.\textit{\ \ \ \ \ \ }$_{\square }\medskip $

\noindent \textbf{Fact 5.5 }\textit{Let }$M$\textit{\ be a lattice
isomorphic to }$N_{5}$\textit{, i.e., }$M=\{0,a,b,c,1\}$\textit{, where }$0$%
\textit{\ and }$1$\textit{\ are its least and greatest elements,
respectively, }$b<c$\textit{, }$a+b=a+c=1$\textit{, and }$a\wedge b=a\wedge
c=0$\textit{. Then, }$End(M)$\textit{\ is not a right CP-semiring.\medskip }

\noindent \textbf{Proof. }Almost verbatim repeating the proof of Fact~5.4,
one only needs to show that the equality $1_{S}=ru+v$, where $1_{S}\in
S:=End(M)$, is impossible when $r\in \{e_{0,c},\,e_{0,1}\}$ and $0\neq u\in
S $. But the latter is almost obvious: Indeed, if $u(a)\neq 0$, then $%
ru(a)\in \{c,1\}$ and, hence, $(ru+v)(a)\in \{c,1\}$ for all $v\in S$, and
therefore, $ru+v\neq 1_{S}$; If $u(a)=0$, then $a+b=1$ implies $u(b)\neq 0$,
$ru(b)\in \{c,1\}$, and, therefore, $(ru+v)(b)\in \{c,1\}$ for any $v\in S$,
and one gets again $ru+v\neq 1_{S}$.\textit{\ \ \ \ \ \ }$_{\square
}\medskip $

\noindent \textbf{Fact 5.6 }\textit{Let }$M$\textit{\ be a chain }$%
1=x_{0}>x_{1}>x_{2}>\ldots >x_{n}>\ldots >0$\textit{. Then, }$S:=End(M)$%
\textit{\ is not a right CP-semiring.\medskip }

\noindent \textbf{Proof. }Clearly, the partition $M=\bigcup_{i=0}^{\infty
}\{x_{2i},x_{2i+1}\}\cup \{0\}$ defines a congruence $\theta $ on $M$. For
every $s,s^{\prime }\in S$, let $N(s,s^{\prime }):=\{s(m)\in M:s(m)\neq
s^{\prime }(m)\}\cup \{s^{\prime }(m)\in M:s(m)\neq s^{\prime }(m)\}$ and $%
\Theta $ be the relation on $S$ defined for all $m\in M$ and $|N(s,s^{\prime
})|<\infty $ as follows: $s\ \Theta $ $s^{\prime }$ $\Longleftrightarrow $ $%
s(m)\ \theta $ $s^{\prime }(m)$. Then, for all $s,s^{\prime },s^{\prime
\prime }\in S$, one may readily verify all inclusions $N(s,s^{\prime \prime
})\subseteq N(s,s^{\prime })\cup N(s^{\prime },s^{\prime \prime })$,
$N(ss^{\prime \prime },s^{\prime }s^{\prime \prime })\subseteq
N(s,s^{\prime })$, and $N(s+s^{\prime \prime },s^{\prime }+s^{\prime \prime
})\subseteq N(s,s^{\prime })$. Whence, keeping them in mind, it is easy to
see that $\Theta $ is a congruence on $S_{S}$.

Now we will show that $\overline{S}:=S/\Theta $ is not a retract of $S_{S}$,
and, hence, $S$ is not a right CP-semiring. Indeed, if it is not a case, by
Lemma~3.7, there exist $e\in I^{\times }(S)$ and an $S$-isomorphism $\psi
\colon \overline{S}\rightarrow eS$ such that $e\ \Theta $ $1$ and $\psi (%
\overline{1})=e$. From the latter, $s\ \Theta $ $s^{\prime }$ $%
\Longleftrightarrow $ $es=es^{\prime }$. Also, one may see right away that $%
x_{i}\ \theta $ $x_{j}$ $\Longleftrightarrow $ $e_{0,x_{i}}\ \Theta $ $%
e_{0,x_{j}}$ and, hence, $x_{i}\ \theta $ $x_{j}$ $\Longleftrightarrow $ $%
ee_{0,x_{i}}=ee_{0,x_{j}}$ $\Longleftrightarrow $ $e(x_{i})=e(x_{j})$.
However, the latter implies $x_{i}\in N(e,1)$ for each $i$; and therefore, $%
N(e,1)$ is infinite, what contradicts to $e\ \Theta $ $1$.\textit{\ \ \ \ \
\ }$_{\square }\medskip $

\noindent \textbf{Proposition 5.7 }\textit{Let }$\theta $\textit{\ be a
congruence on }$M$\textit{\ and }$\overline{M}:=M/\theta $\textit{. If }$%
S:=End(M)$\textit{\ is a right CP-semiring, then so is }$\overline{S}:=End(%
\overline{M})$\textit{.}\medskip

\noindent \textbf{Proof. }So, let $S$ be a right CP-semiring. First notice
that the congruence $\theta $ naturally produces the congruence $\theta F$
on $S_{S}$ defined as follows:
\begin{equation*}
s\ \theta F\text{ }s^{\prime }\Longleftrightarrow \forall m\in M:\text{ }%
s(m)\ \theta \text{ }s^{\prime }(m)\text{.}
\end{equation*}%
Also, $\forall x,y\in M:$ $x\ \theta $ $y$ $\Leftrightarrow $ $e_{0,x}\
\theta F$ $e_{0,y}$.

By Lemma~3.7, there exists an element $e\in I^{\times }(S)$ such that $%
S/(\theta F)\cong eS$ and $1$ $\theta F$ $e$. By Proposition~5.1, the
semiring $R=eSe$ is a right CP-semiring too. So, to complete the proof, it
is enough to find a semiring isomorphism $\gamma \colon R\rightarrow
\overline{S}$.

Let $\alpha \colon M\rightarrow \overline{M}$ be the canonical surjection,
and $(\gamma (r))(\overline{x}):=\alpha (r(x))$ for every $r\in R$ and $%
\overline{x}\in \overline{M}$. We shall show that $\gamma $ actually defines
a desired semiring isomorphism $\gamma \colon R\rightarrow \overline{S}$.

First of all, as in the proof of Fact~5.6, one may easily see that the
following is true: $\forall x,y\in M:$ $x\ \theta $ $y$ $\Leftrightarrow $ $%
e(x)=e(y)$. Hence, for every $r\in R$, the mapping $\gamma (r)\colon
\overline{M}\rightarrow \overline{M}$ is well-defined. For $r$ and $\alpha $
are homomorphisms, $\gamma (r)$ is an endomorphism, \textit{i.e.}, $\gamma
(r)\in \overline{S}$.

Furthermore, for any $\overline{x}\in \overline{M}$ and $r,r^{\prime }\in R$%
, we have $[\gamma (r+r^{\prime })](\overline{x})=[\alpha (r+r^{\prime
})](x)=\alpha r(x)+\alpha r^{\prime }(x)=[\gamma (r)+\gamma (r^{\prime })](%
\overline{x})$, and hence, $\gamma (r+r^{\prime })=\gamma (r)+\gamma
(r^{\prime })$; also, we have $\gamma (rr^{\prime })(\overline{x})=\alpha
(rr^{\prime })(x)=\gamma (r)(\overline{r^{\prime }(x)})=\gamma (r)(\alpha
(r^{\prime }(x)))=\gamma (r)(\gamma (r^{\prime })(\overline{x}))=(\gamma
(r)\gamma (r^{\prime }))(\overline{x})$, and hence, $\gamma (rr^{\prime
})=\gamma (r)\gamma (r^{\prime })$; obviously, $\gamma (0)=0_{\overline{S}}$%
. So, $\gamma $ is a semiring homomorphism.

Let $r,r^{\prime }\in R$ and $r\neq r^{\prime }$. Then, for some $x\in M$,
we have $r(x)\neq r^{\prime }(x)$. Whence, $(r(x),r^{\prime }(x))\notin
\,\theta $ (otherwise, $r(x)=e(r(x))=e(r^{\prime }(x))=r^{\prime }(x)$); and
therefore, $\gamma (r)(\overline{x})=\alpha (r(x))\neq \alpha (r(x))=\gamma
(r^{\prime })(\overline{x})$ and $\gamma $ is an injection.

Finally, let $\overline{s}\in \overline{S}$ and $\tilde{e}$ be the
homomorphism from $\overline{M}$ to $M$ such that $\tilde{e}(\overline{x}%
):=e(x)$. In particular, $e\tilde{e}=\tilde{e}$; moreover, since $x\ \theta $
$e(x)$ for all $x\in M$, we also have $\alpha e=\alpha $ and $\alpha \tilde{e%
}=1_{\overline{M}}$. Then, putting $r:=\tilde{e}\overline{s}\alpha $, we
have $r=ere\in R$ and
\begin{equation*}
\gamma (r)(\overline{x})=\alpha r(x)=\alpha \tilde{e}\overline{s}\alpha (x)=%
\overline{s}(\overline{x})
\end{equation*}%
for every $\overline{x}\in \overline{M}$, hence, $\gamma (r)=\overline{s}$
and $\gamma $ is a surjection as well.\textit{\ \ \ \ \ \ }$_{\square
}\medskip $

\noindent \textbf{Proposition 5.8 }\textit{\ If }$S=End(M)$\textit{\ is a
right CP-semiring, then} $M$ \textit{is a finite distributive
lattice.\medskip }

\noindent \textbf{Proof. }By Lemma~5.3, $M$ is a bounded lattice. Let $%
L\subseteq M$ be a sublattice and $U(x):=\{a\in L\,|\,x\leq a\}$ for every $%
x\in M$. Clearly, $U(x)$ is either an empty set, or a sublattice in $L$. In
the latter, let $z_{x}$ denote its least element if it exists.

First suppose that $L$ is finite and $\phi \colon M\rightarrow L$ be the
map, defined for every $x\in M$ as follows:
\begin{equation*}
\phi (x):=\{_{1_{L}\text{ otherwise \ }}^{z_{x}\text{ if }U(x)\neq \emptyset
}\text{.}
\end{equation*}%
One can readily verify that $\phi $ is a surjective homomorphism and it
induces the natural congruence $\sim _{\phi }$ on $M$ defined for all $%
x,y\in M$ by formula:
\begin{equation*}
x\sim _{\phi }y\Longleftrightarrow \phi (x)=\phi (y)\text{,}
\end{equation*}%
and $M/{\sim _{\phi }}\cong L$. From this observation, Proposition~5.7, and
Facts~5.4 and 5.5, we conclude that $M$ has no sublattices isomorphic to $%
M_{3}$ or $N_{5}$, and therefore, by \cite[Theorem IX.2]{birkhoff:lathe} or
\cite[Theorem 2.1.1]{gratzer:glt}, $M$ is a distributive lattice.

Now let $L$ be a chain $a_{0}>a_{1}>a_{2}>\dots >a_{n}>\ldots $ Without loss
of generality, we may assume that $a_{0}$ is the greatest element of $M$.
Then, $U(x)=L$ or $U(x)$ has the least element for every $x\in M$. Let $%
L^{\prime }$ denote the chain from Fact~5.6 and define the surjective
homomorphism $\phi \colon M\rightarrow L^{\prime }$ as follows: $\phi
(x):=x_{k}$ when $a_{k}$ is the least element of $U(x)$, and $\phi
(x):=0_{L^{\prime }}$ for $U(x)=L$. It induces the natural congruence $\sim
_{\phi }$ on $M$ such that $M/{\sim _{\phi }}\cong L^{\prime }$. From the
latter, Proposition~5.7 and Fact~5.6, we get that $M$ satisfies the
descending chain condition. Using this fact and Lemma 5.2, we get that $M$
has a finite maximal subchain, that is, $M$ is a distributive lattice of a
finite length, and therefore (see, for example, \cite[Exercises 6, Page 127]%
{skor:eolt}), $M$ is finite.\textit{\ \ \ \ \ \ }$_{\square }\medskip $

For the reader's convenience, remind some fundamental concepts and facts
regarding the Morita equivalence of semirings from \cite{kat:thcos} and \cite%
{kn:meahcos} that we will use in sequence. \ Thus, two semirings $T$ and $S$
are said to be \textit{Morita equivalent} if the semimodule categories $_{T}%
\mathcal{M}$ and $_{S}\mathcal{M}$ are equivalent categories, \textit{i.e.},
there exist two additive functors $F:$ $_{T}\mathcal{M}$ $\longrightarrow $ $%
_{S}\mathcal{M}$ and $G:$ $_{S}\mathcal{M}\longrightarrow $ $_{T}\mathcal{M}$
and natural isomorphisms $\eta :GF\longrightarrow Id_{_{T}\mathcal{M}}$ and $%
\xi :FG\longrightarrow Id_{_{S}\mathcal{M}}$. Two semirings $T$ and $S$ are
Morita equivalent iff the semimodule categories $\mathcal{M}_{T}$ and $%
\mathcal{M}_{S}$ are equivalent categories \cite[Theorem 4.12]{kn:meahcos}.
A left semimodule $_{T}P\in |_{T}\mathcal{M}|$ is a \textit{generator} in
the category $_{T}\mathcal{M}$ if the regular semimodule $_{T}T\in |_{T}%
\mathcal{M}|$ is a retract of a finite direct sum $\oplus _{i}P$ of the
semimodule $_{T}P$; and a left semimodule $_{T}P\in |_{T}\mathcal{M}|$ is a
\textit{progenerator} in $_{T}\mathcal{M}$ if it is a finitely generated
projective generator. Finally, two semirings $T$ and $S$ are Morita
equivalent iff there exists a progenerator $_{T}P\in |_{T}\mathcal{M}|$ in $%
_{T}\mathcal{M}$ such that the semirings $S$ and $\mathrm{End}$$(_{T}P)$
are isomorphic \cite[Theorem 4.12]{kn:meahcos}.

Now let $M$ be a finite distributive lattice. Then, the set
\begin{equation*}
T(M):=\{m\in M\,|\,M=m^{\triangledown }\cup m^{\vartriangle }\}
\end{equation*}%
is obviously a chain containing $0$ and $1$. Also, let $[a,b]:=\{x\in
M\,|\,a\leq x\leq b\}$ denote the intervals defined for all $a,b\in M$ and $%
a\leq b$. We say that an interval $[a,b]$ is \textit{simple} if $%
[a,b]=\{a,b\}$. Using these notions, our next theorem provides us with a
full description of finite distributive lattices whose endomorphism
semirings are right CP-semirings.\medskip

\noindent \textbf{Theorem 5.9 }\textit{For a finite distributive lattice }$M$%
\textit{, the following conditions are equivalent:}

\textit{(1) }$End(M)$\textit{\ is a right CP-semiring;}

\textit{(2)}$\ M/\theta $\textit{\ is a distributive lattice for any }$%
\theta \in Cong(M)$\textit{;}

\textit{(3) If }$t,t^{\prime }\in T(M)$\textit{\ and }$[t,t^{\prime
}]|_{T(M)}$\textit{\ is simple, then }$[t,t^{\prime }]$\textit{\ is either
simple, or isomorphic to }$\mathbf{B}^{2}$\textit{.}\medskip

\noindent \textbf{Proof. }(1) $\Longrightarrow $ (2). This implication
follows immediately from Propositions~5.7 and~5.8.\smallskip

(2) $\Longrightarrow $ (1). Let $S:=End(M)$. Since $M$ is a finite
distributive lattice, from \cite[Theorem 5.3]{hk:tcos} (or, see, \cite[Fact
5.9]{kat:thcos}) it follows that the semimodule $_{\mathbf{B}}M$ is a
progenerator in the category $_{\mathbf{B}}\mathcal{M}$. Also, by \cite[%
Corollary 3.13]{kn:meahcos}, the $\mathbf{B}$-semimodule $M^{\ast }:=Hom_{%
\mathbf{B}}(M,\mathbf{B})$ is a progenerator in the category $_{\mathbf{B}}%
\mathcal{M}$, and, using \cite[Proposition 3.12]{kn:meahcos}, one gets $_{%
\mathbf{B}}M\cong $ $_{\mathbf{B}}(M^{\ast })^{\ast }:=Hom_{\mathbf{B}%
}(M^{\ast },\mathbf{B})$ and $S\cong End(_{\mathbf{B}}M^{\ast })$; that is, $%
S$ is Morita equivalent to $\mathbf{B}$ via the progenerator $_{\mathbf{B}%
}M^{\ast }$ in $_{\mathbf{B}}\mathcal{M}$. Then, as showed in the proof of
\cite[Theorem 4.12]{kn:meahcos}, the functors $F:$ $\mathcal{M}%
_{S}\rightleftarrows $ $\mathcal{M}_{\mathbf{B}}:G$, defined by $%
F(A)=A\otimes _{S}(M^{\ast })^{\ast }\cong A\otimes _{S}M$ and $%
G(B)=B\otimes _{\mathbf{B}}M^{\ast }$, establish an equivalence between the
semimodule categories.

Let $A$ be a cyclic right $S$-semimodule. Then, there exists a surjective $S$%
-homomorphism $f\colon S\rightarrow A$. By \cite[Lemma 4.7]{kn:meahcos}, the
homomorphism $F(f):M\cong F(S)=S\otimes _{S}M\longrightarrow F(A)$ is
surjective in $\mathcal{M}_{\mathbf{B}}$ as well, and hence, there exists a
congruence $\theta $ on $M$ such that $F(A)\cong M/\theta $. By hypothesis
(2), $F(A)$ is a finite distributive lattice; and hence, by~\cite[Theorem~5.3%
]{hk:tcos} (or, see, \cite[Fact 5.9]{kat:thcos}), it is a projective $%
\mathbf{B}$-semimodule. Applying \cite[Lemma~4.10]{kn:meahcos}, we obtain
that $A\cong G(F(A))$ is a projective right $S$-semimodule, and therefore, $%
S $ is a right CP-semiring. \smallskip

(2) $\Longrightarrow $ (3). Let $t,t^{\prime }\in T(M)$ and $[t,t^{\prime
}]|_{T(M)}$ is simple. The mapping $\alpha \colon M\rightarrow \lbrack
t,t^{\prime }]$ such that $\alpha (x):=t$ for $x\in t^{\triangledown }$, $%
\alpha (x):=x$ for $x\in \lbrack t,t^{\prime }]$, and $\alpha (x):=t^{\prime
}$ otherwise, obviously is a surjective homomorphism. Whence, $[t,t^{\prime
}]\cong M/{\sim _{\alpha }}$ for the natural congruence $\sim _{\alpha }$ on $%
M$. It is clear that, if $M$ satisfies (2), then every quotient semimodule
of $M$ satisfies (2) as well; so, without loss of generality, we may assume
that $M=[t,t^{\prime }]$ and, hence, $T(M)=\{0,1\}$. We shall show that $%
M\cong \mathbf{B}$, or $M\cong \mathbf{B}^{2}$.

For $1<|M|$ $<\infty $, $M$ has, at least, one atom --- a minimal element in
$M\backslash \{0\}$. On the other hand, it is easy to see that $M$ has, at
most, two atoms: Indeed, if $a$, $b$ and $c$ are different atoms of $M$,
then, there exists the surjective homomorphism $\phi :M\longrightarrow M_{3}$
defined as follows: $\phi (0):=0$, $\phi (a):=a$, $\phi (b):=b$, $\phi
(c):=c $, and $\phi (x):=1$ otherwise, that, since $M_{3}$ is not a
distributive lattice, contradicts (2). Thus, $M$ has one or two atoms.

Let $a\in M$ be the unique atom. Then, $a\in T(M)$, and therefore, $a=1$ and
$M=\{0,1\}\cong \mathbf{B}$.

Now consider the remaining case when $M$ has precisely two atoms --- $a$ and
$b$. Then, we have $a\leq c$, or $b\leq c$ for each nonzero $c\in M$ and,
therefore, $[0,a+b]=\{0,a,b,a+b\}$ (since otherwise, $\{0,a,b,c,a+b\}%
\subseteq M$ is a sublattice which is isomorphic to $N_{5}$ and $M$ would
not be distributive), and consider the following two cases.

If $a+b\in T(M)$, then $a+b=1$ and $M=[0,1]=[0,a+b]\cong \mathbf{B}^{2}$.

Finally, suppose $a+b\not\in T(M)$. Then, there exists $c\in M$ such that $%
a+b$ and $c$ are incomparable with each other. Hence, just one of the
inequalities $a<c$ and $b<c$ is true, say, for instance, the second one.
Without loss of generality, we may also assume that $[b,c]$ is simple.
However, one can easily see that in this case there exists the surjective
homomorphism $\phi :M\longrightarrow N_{5}$ defined as follows: $\phi (0):=0$%
, $\phi (a):=a$, $\phi (b):=b$, $\phi (c):=c$, and $\phi (x):=1$ otherwise,
that, since $N_{5}$ is not a distributive lattice, contradicts (2).

(3) $\Longrightarrow $ (2). Obviously, if a lattice satisfies (3), then its
every quotient $\mathbf{B}$-semimodule satisfies (3), too. Also, each
lattice satisfying (3) is distributive for it has no sublattices isomorphic
to $M_{3}$ or $N_{5}$. From these observations, we end the proof.\textit{\ \
\ \ \ \ }$_{\square }\medskip $

Our next result, significantly extending Theorem 5.9, gives a full
description of semimodules over arbitrary Boolean algebras whose
endomorphism semirings are right CP-semirings.\medskip

\noindent \textbf{Theorem 5.10 }\textit{Let }$M$\textit{\ be a nonzero
semimodule over an arbitrary Boolean algebra }$B,$ $S:=End(M_{B})$%
\textit{\ and }$J:=Ann_{B}(M)$\textit{. Then, }$S$\textit{\ is a CP-semiring
iff }$B/J$\textit{\ is a finite Boolean algebra and }$Ma$\textit{\ for every
atom }$a\in B/J$ \textit{is a finite distributive lattice satisfying the
equivalent conditions of Theorem 5.9.\medskip }

\noindent \textbf{Proof. }$\Longrightarrow $. Let $S$ be a right
CP-semiring, and $\overline{B}:=B/J$. One sees right away that $M$ is also a
right $\overline{B}$-semimodule and $S\cong End(M_{\overline{B}})=:\overline{%
S}$. For $S$ is a right CP-semiring, $\overline{S\text{ }}$is a right
CP-semiring as well. It is easy to see that $\overline{B}$ is an Boolean
algebra, too, and we will show that $|\overline{B}|$ $<\infty $.

Indeed, if $\overline{B}$ is infinite, then, by \cite[Theorem 32]{Sal1988},
it contains a countable set of orthogonal idempotents $\{e_{1},e_{2},...\}$.
For each $e_{i}$, let $\alpha _{i}:M\longrightarrow M$ be the $\overline{B}$%
-homomorphism given by the formula $\alpha _{i}(m)=me_{i}$. It is clear that
each $\alpha _{i}$ is a nonzero idempotent element in $\overline{S}$ and $%
\alpha _{i}\alpha _{j}=0$ for all $i,j$ such that $j\neq i$. On the other
hand, for $\overline{S}$ is a right CP-semiring and by Lemma~3.9, all but
finite number of elements of the family $\{\alpha _{i}\}_{i\in \boldsymbol{N}%
}$ should be of zero. Thus, $\overline{B}$ is a finite Boolean algebra.

Let $\{e_{1},e_{2},...,e_{n}\}$ be the set of all atoms of $\overline{B}$.
One sees right away that $e_{i}\overline{B}=\{0,e_{i}\}\cong \mathbf{B}$ for
each $i$, and $\overline{B}=\overline{B}e_{1}\oplus \overline{B}e_{2}\oplus
...\oplus \overline{B}e_{n}$ and $M=Me_{1}\oplus Me_{2}\oplus ...\oplus
Me_{n}$; moreover, each $Me_{i}$ is, in fact, a right $e_{i}\overline{B}$%
-semimodule. Whence,
\begin{equation*}
\overline{S}\cong End(Me_{1})\oplus End(Me_{2})\oplus ...\oplus End(Me_{n})%
\text{.}
\end{equation*}%
Then, using Corollary 3.6, Proposition~5.8 and Theorem~5.9, one concludes
this implication.

$\Longleftarrow $. Taking into consideration that $S\cong End(Me_{1})\oplus
...\oplus End(Me_{n})$, where each $e_{i}$ is an atom of $\overline{B}$, and
applying Corollary 3.6, Proposition~5.8 and Theorem~5.9, we end the proof.%
\textit{\ \ \ \ \ \ }$_{\square }\medskip $

From Theorems~4.4, 4.16, and \cite[Theorem 5.10]{knz:ososacs}, it is easy to
see that for a semiring to be `ideal-simple' and a `CP-semiring' in general
are \textquotedblleft independent\textquotedblright\ properties --- there
exist both CP-semirings that are ideal-simple and that are not ideal-simple.
In light of this, the following applications of Theorem 5.9, completely
describing ideal-simple CP-semirings, are certainly of the interest and
conclude a list of the central results of the paper.\medskip

\noindent \textbf{Theorem 5.11 } \textit{For a semiring }$S$\textit{, the
following conditions are equivalent:}

\textit{(1) }$S$\textit{\ is an ideal-simple right CP-semiring;}

\textit{(2) }$S$\textit{\ is a simple right CP-semiring;}

\textit{(3) }$S\cong M_{n}(D)$\textit{\ for some division ring }$D$\textit{,
or }$S\cong $\textit{\ }$End(M)$\textit{, where }$M$\textit{\ is a finite
distributive lattice satisfying the equivalent conditions of Theorem
5.9.\medskip }

\textbf{Proof. }(1) $\Longrightarrow $ (2). Assume that $S$ is an
ideal-simple right CP-semiring. By Theorem~3.11, $S$ is a semisimple ring,
or $S$ is a zeroic right CP-semiring. In the first case, the implication is
trivial. So, let $S$ be a zeroic right CP-semiring. By Proposition~3.10, $S$
is additively $\pi $-regular, \textit{i.e.}, there exist a natural number $%
n\geq 1$ and an element $x\in S$ such that $n1=n1+x+n1$. For $S$ is a
zerosumfree semiring, $(n1+x)\in S$ is a nonzero additively idempotent
element and $I^{+}(S)$ is a nonzero ideal of $S$. Whence, $S=I^{+}(S)$ and,
by \cite[Proposition 3.4]{gw:oeics} (or \cite[Proposition 23.5]{golan:sata}%
), $S$ can be embedded in some complete semiring. From the latter and \cite[%
Proposition 22.27]{golan:sata}, it follows that $S$ can be considered as a
subsemiring of a semiring $T$ with an infinite element $\infty $. For $S$ is
a right CP-semiring and Lemma~3.7, there exist $e\in I^{\times }(S)$ and an $%
S$-isomorphism $\varphi :\infty S\longrightarrow eS$, with $e=\varphi
(\infty )$, and, by Proposition~5.1, $eSe$ is a right CP-semiring, too.
Moreover, for each $s\in S$, $ese+e=(es+e)e=(\varphi (\infty )s+\varphi
(\infty ))e=(\varphi (\infty s)+\varphi (\infty ))e=\varphi (\infty s+\infty
)e=\varphi (\infty )e=e^{2}=e$, and hence, $eSe$ is a Gelfand semiring. And
therefore, by Theorem~4.9, $eSe$ is a finite Boolean algebra. On the other
hand, by \cite[Proposition 5.3]{knz:ososacs}, $eSe$ is an ideal-simple
semiring as $S$ is ideal-simple; and hence, $eSe\cong \mathbf{B}$. Also, for
$S$ is an ideal-simple semiring, we have that $SeS=S$; and therefore, by
\cite[Proposition 5.2]{knz:ososacs}, $S$ is Morita equivalent to $\mathbf{B}$
and, by \cite[Theorem 5.6]{knz:ososacs}, $S$ is a simple semiring.

(2) $\Longrightarrow $ (3). Let $S$ be a simple right CP-semiring. By
Theorem~3.11, $S$ is a semisimple ring, or $S$ is a zeroic right
CP-semiring. In the first case, it is trivial that $S\cong M_{n}(D)$ for
some division ring $D$. So, let $S$ be a zeroic right CP-semiring. Let $%
S^{\ast }:=Hom_{\mathbb{N}}(S,\mathbf{B})$, then $S^{\ast }$ is a right $S$%
-semimodule with the scalar multiplication given as: $\phi s(x):=\phi (sx)$
for all $x,s\in S$ and $\phi \in S^{\ast }$. Obviously, $(S^{\ast },+,0)$ is
an upper semilattice. For $S$ is also a zerosumfree semiring, the $\mathbb{N}
$-homomorphism $\varphi _{0}:S\longrightarrow \mathbf{B}$, given by $\varphi
_{0}(0):=0$, and $\varphi _{0}(x):=1$ for all $0\neq x\in S$, is the
greatest element in $S^{\ast }$; and let us consider the cyclic right $S$%
-semimodule $\varphi _{0}S$. By Zorn's lemma, there exists a maximal
congruence $\sim $ on the $S$-semimodule $\varphi _{0}S$, and let $%
M:=\varphi _{0}S/\sim $. Obviously, the quotient semimodule $M$ has only the
trivial congruences. Let $N$ be a subsemimodule of $M$; then $M$ has only
the trivial congruences and the Bourne congruence $\equiv _{N}$ is either
the diagonal congruence, or the universal one. If $\equiv _{N}$ is the
diagonal congruence, then $x\equiv _{N}0$ for any $x\in N$ and, hence, $N=0$%
. Otherwise, $m\equiv _{N}m^{\prime }$ for all $m,m^{\prime }\in M$; in
particular, we have $\overline{\varphi }_{0}\equiv _{N}0$, \textit{i.e.},
there exist elements $x,y\in N$ such that $\overline{\varphi }_{0}=\overline{%
\varphi }_{0}+x=0+y=y\in N$, which implies that $N=M$. Therefore, $M$ is a
minimal right $S$-semimodule. Then, for $S$ is a right CP-semiring and
Lemma~3.7, there exists $e\in I^{\times }(S)$ such that $eS\cong M$ as right
$S$-semimodules. So, $eS$ is a projective minimal right ideal of $S$.
Whence, by \cite[Theorem 5.10]{knz:ososacs}, we get that $S\cong End(M)$ for
some finite distributive lattice $M$, and $M$ satisfies the equivalent
conditions of Theorem~5.9.

(3) $\Longrightarrow $ (1). This implication follows immediately from
Theorems~4.4 and~5.9, and \cite[Theorem 5.10]{knz:ososacs}.\textit{\ \ \ \ \
\ }$_{\square }\medskip $

As shown in \cite[Examples 3.8(b)]{zumbr:cofcsswz}, the concepts of
`congruence-simpleness' and `ideal-simpleness' for finite semirings are not
the same. However, it is not a case for CP-semirings:\medskip

\noindent \textbf{Corollary 5.12 }\textit{For a finite semiring }$S,$\textit{%
\ the following conditions are equivalent:}

\textit{(1) }$S$\textit{\ is an ideal-simple right CP-semiring;}

\textit{(2) }$S$\textit{\ is a congruence-simple right CP-semiring;}

\textit{(3) }$S\cong M_{n}(F)$\textit{\ for some finite field }$F$\textit{,
or }$S\cong End(M)$\textit{, where }$M$\textit{\ is a finite distributive
lattice satisfying the equivalent conditions of Theorem~5.9.\medskip }

\noindent \textbf{Proof. }The implications (1) $\Longrightarrow $ (2) and
(3) $\Longrightarrow $ (1) follow from Theorem~5.11.

(2) $\Longrightarrow $ (3). By Theorem~3.11, $S$ is a finite semisimple
ring, or a finite zeroic right CP-semiring. For $S$ is congruence-simple, in
the first \textquotedblleft scenario,\textquotedblright\ it is clear that $%
S\cong M_{n}(F)$ for some finite field $F$.

Thus, consider the case when $S$ is a finite zeroic right CP-semiring. By
\cite[Theorem 1.7]{zumbr:cofcsswz}, there exists a nonzero finite $\mathbf{B}
$-semimodule $M$ such that $S$ is a subsemiring of $End(M)$ containing all
endomorphisms $e_{x,y}$, $x,y\in M$; and we will show that $M$ is a finite
distributive lattice.

Indeed, suppose that $L$ is a sublattice of $M$ isomorphic either to the
lattices $M_{3}$ or $N_{5}$. Then, as was done in the proof of Proposition
5.8, there are the surjective homomorphism $\phi :M\longrightarrow L$ and
the natural congruence $\theta $ on $M$ such that $\overline{M}:=M/\theta
\cong L$. And we have (see also Proposition 5.7) the congruence $\theta F$
on $S_{S}$ defined as follows:
\begin{equation*}
s\ \theta F\text{ }s^{\prime }\Longleftrightarrow \forall m\in M:\text{ }%
s(m)\ \theta \text{ }s^{\prime }(m)\text{.}
\end{equation*}%
For $S$ is a right CP-semiring and Lemma~3.7, there exists an element $e\in
I^{\times }(S)$ such that $S/(\theta F)\cong eS$ and $1$ $\theta F$ $e$,
and, by Proposition~5.1, the semiring $R:=eSe$ is a right CP-semiring, too.

Analogously as was done in the proof of Proposition~5.7, one may readily
verify that the following is true: $\forall\, x,y\in M:$ $x\ \theta $ $y$ $%
\Leftrightarrow $ $e(x)=e(y)$; $\forall\, x,y\in M:$ $x$ $\theta $ $e(x)$; the
mapping $\gamma :R\longrightarrow End(\overline{M})$, defined as follows: $%
\gamma (r)(\overline{x}):=\alpha r(x)$ for all $r\in R$ and $\overline{x}\in
\overline{M}$, where $\alpha :M\longrightarrow \overline{M}$ is the
canonical projection, is an injective homomorphism of semirings. Then, the
semiring $\overline{S}:=\gamma (R)$ is obviously isomorphic to $R$ and,
hence, $\overline{S}$ is a right CP-semiring, too.

Notice that $\overline{m}=$ $\overline{0}$ iff $e(m)=e(0)=0$ for all $m\in M$%
. Whence, $\gamma (ee_{0,x}e)=e_{\overline{0},\overline{x}}$, and the
semiring $\overline{S}$ contains all endomorphisms $e_{\overline{0},%
\overline{x}}$ for all $\overline{x}\in \overline{M}$. Now, analogously as
was done in the proofs of Fact 5.4, or Fact 5.5, respectively, we have that $%
\overline{S}$ is not a right CP-semiring. Thus, $M$ is a finite distributive
lattice, and, using \cite[Proposition 4.9 and Remark 4.10]{zumbr:cofcsswz},
we obtain that $S=End(M)$ and, by Theorem~5.11, end the proof.\textit{\ \ \
\ \ \ }$_{\square }\medskip $

\section{Some conclusive remarks and problems}

As was briefly mentioned in Section 3, studying of representations of
semimodules over a semiring $S$ as colimits of full c-diagrams of the
regular semimodules $S_{S}$, in our opinion, constitutes a very interesting
and promising and, we believe, innovative area of research --- homological
structure theory of semirings --- having, of course, the \textquotedblleft
bloody\textquotedblright\ connections, as was mentioned earlier too, with
the homological characterization of semirings as well as with numerous
diverse areas of modern algebra, topology and model theory. Therefore, we
would be glad to motivate and encourage our potential readers to join us and
to further the homological structure theory (of semirings) significantly
wider in the following two ways: First, to consider colimits of different
types/sorts of diagrams (not only full c-diagrams) whose objects, in turn,
are from very well established specified classes of (semi)modules; Second,
to develop the homological structure theory in some nonadditive important
settings --- first of all, such as $S$-acts (see \cite{kilp-kn-mik:maac})
and Grothendieck toposes of (pre)sheaves (see \cite{maclmoer:sigal}, or \cite%
{kashschap:cas}).

Finally, a few, in our view, interesting specific problems closely related
to the considerations in this paper.

As was shown in \cite[Proposition 4.1]{aikn:ovsasaowcsai}, every zerosumfree
CI-semiring possesses an infinite element. As to zerosumfree CP-semirings,
by Proposition 3.10, they are zeroic; however, all our examples of such
semirings are semirings with infinite elements, and we have the following
conjecture\medskip\

\noindent \textbf{Conjecture 6.1. }Every zerosumfree CP-semiring possesses
an infinite element.\medskip

In contrast with the case of rings, the classes of CI-semirings and
CP-semirings are different ones (compare, for example, Theorem 4.16 and \cite%
[Theorem 4.20]{aikn:ovsasaowcsai}). Therefore, it would be interesting to
see how far the classical result --- Theorem 3.3 --- of the ring case can be
extended for semirings, namely:\medskip

\noindent \textbf{Problem 6.2.} Describe semirings $S$ such that full
c-diagrams of every right (left) semimodule over which are always injective and
projective ones.\medskip

Also, it would be interesting to see how far the \textquotedblleft
symmetrical\textquotedblright\ part of Theorem 3.3 can be extended for
semirings, namely:\medskip

\noindent \textbf{Problem 6.3.} Describe semirings $S$ which are right
CP-semirings (CI-semirings) iff they are left CP-semirings
(CI-semirings).\medskip

Finally, describing CP-semirings (CI-semirings)\ of endomorphisms of
semimodules over various ground semirings $S$ in the spirit of Section 5, in
our view, constitutes a quite interesting, promising, and not trivial at
all, direction for the furthering research initiated in Section 5, for
example:\medskip\

\noindent \textbf{Problem 6.4.} Describe semimodules $M$ over a distributive
lattice $D$ whose semiring of endomorphisms $End(M_{D})$ are CP-semirings.

\end{document}